\documentclass{article}

\usepackage{arxiv}
\usepackage[T1]{fontenc}
\usepackage{amsmath,amssymb,amsthm,mathtools,bm,mathrsfs}
\usepackage{booktabs}
\usepackage{enumitem}
\usepackage{mdframed}
\usepackage{microtype}
\usepackage[square,numbers]{natbib}
\usepackage{xcolor}
\usepackage[colorlinks=true,linkcolor=blue,citecolor=blue,urlcolor=blue]{hyperref}

\usepackage[nameinlink,capitalize]{cleveref}
\crefname{equation}{}{}
\Crefname{equation}{}{}

\allowdisplaybreaks

\newtheorem{theorem}{Theorem}[section]
\newtheorem{proposition}[theorem]{Proposition}
\newtheorem{lemma}[theorem]{Lemma}
\newtheorem{corollary}[theorem]{Corollary}
\newtheorem{remark}[theorem]{Remark}
\newtheorem{example}[theorem]{Example}
\theoremstyle{definition}
\newtheorem{definition}[theorem]{Definition}

\newmdenv[
  linewidth=0.7pt,
  roundcorner=2pt,
  innerleftmargin=8pt,
  innerrightmargin=8pt,
  innertopmargin=2pt,
  innerbottommargin=3pt,
  skipabove=8pt,
  skipbelow=8pt
]{mainresult}

\newcommand{\R}{\mathbb R}
\newcommand{\Sph}{\mathbb S}

\newcommand{\Rn}{\mathbb{R}^n}
\newcommand{\Rexx}{[-\infty,\infty]}
\newcommand{\Rex}{(-\infty,\infty]}
\newcommand{\bB}{\mathbb{B}}
\newcommand{\N}{\mathbb{N}}
\newcommand{\sL}{\mathbb{L}}
\newcommand{\sphL}{\mathbb{S}_{\!\;\mathbb L}}

\newcommand{\iprod}[2]{\langle #1, #2 \rangle}
\newcommand{\norm}[1]{\lVert #1\rVert}

\newcommand{\dist}{\operatorname{dist}}
\newcommand{\dom}{\operatorname{dom}}
\newcommand{\gph}{\operatorname{gph}}
\newcommand{\epi}{\operatorname{epi}}
\newcommand{\ri}{\operatorname{ri}}
\newcommand{\lin}{\operatorname{lin}}
\newcommand{\cl}{\operatorname{cl}}
\newcommand{\Span}{\operatorname{span}}

\newcommand{\co}{\operatorname{conv}}
\newcommand{\rd}{\mathrm{d}}
\newcommand{\degree}{\operatorname{deg}}

\newcommand{\argmin}{\operatorname*{argmin}}

\newcommand{\prox}{\operatorname{prox}}
\newcommand{\env}[2]{\operatorname{env}_{#1}^{#2}}
\newcommand{\sgn}{\operatorname{sgn}}

\newcommand{\proxh}{\mathrm{prox}_{\lambda h}}
\newcommand{\envh}{\env{h}{\lambda}}
\newcommand{\proxth}{\mathrm{prox}_{\lambda \tilde h}}
\newcommand{\envth}{\env{\tilde h}{\lambda}}

\newcommand{\Id}{\operatorname{Id}}

\newcommand{\nminus}{n_-}

\newcommand{\half}{\tfrac12}
\newcommand{\veps}{\varepsilon}

\newcommand{\linq}{{\mathscr L}_q}
\newcommand{\eq}{{\mathscr E}}
\newcommand{\negq}{{\mathscr N}_q}

\title{From Metric to Strong Metric Regularity: A Second Subderivative Criterion for Definable Prox-regular Subdifferentials}
\author{%
Wenqing Ouyang\\
Department of Mathematics and Statistics\\
Toulouse School of Economics\\
1, Esplanade de l'Universit\'e\\
31080 Toulouse Cedex 06, France\\
\texttt{wenqing.ouyang@tse-fr.eu}
\And
Andre Milzarek\thanks{Andre Milzarek was partly supported by the National Natural Science Foundation of China under Grant No.\ W2532005 and by the Guangdong Provincial Key Laboratory of Mathematical Foundations for Artificial Intelligence (2023B1212010001).}\\
School of Data Science\\
The Chinese University of Hong Kong, Shenzhen\\
Guangdong 518172, P.R. China\\
\texttt{andremilzarek@cuhk.edu.cn}}
\date{}

\renewcommand{\headeright}{Preprint}
\renewcommand{\undertitle}{Preprint}
\renewcommand{\shorttitle}{Metric and Strong Metric Regularity of Definable Prox-regular Subdifferentials}
\hypersetup{
  pdftitle={Metric Regularity versus Strong Metric Regularity for Definable Prox-Regular Subdifferentials},
  pdfauthor={Wenqing Ouyang and Andre Milzarek}
}

\begin{document}
\maketitle

\begin{abstract}
We investigate when metric regularity of a subdifferential mapping rules out multiple local inverse branches and therefore upgrades to strong metric regularity. Let $h:\R^n\to(-\infty,\infty]$ be proper, lower semicontinuous, definable, prox-regular, and subdifferentially continuous at $(\bar x,\bar v)\in\gph(\partial h)$. We characterize this ``upgrade'' through the second subderivative $q = \mathrm{d}^2h(\bar x|\bar v)$ of $h$. If $\partial h$ is metrically regular, then it is strongly metrically regular if and only if the set of negative curvature directions of $q$ on the unit sphere of $\Span(\dom(q))$ has Euler characteristic $0$ or $2$. We also prove a simpler and directly verifiable sufficient condition: nonnegativity of $q$ on the lineality space $\dom(q)\cap(-\dom(q))$ forces a quadratic lower bound for $q$ on its entire domain and hence yields strong metric regularity. In particular, if $\bar x$ is a local minimizer and $\bar v=0$, the second subderivative is nonnegative on its domain $\dom(q)$, which proves a conjecture of Drusvyatskiy, Mordukhovich, and Nghia in the definable setting. The derived principle strengthens this conjecture by replacing local minimality with a weaker second-order condition on the lineality space.
\end{abstract}

\renewcommand{\keywordname}{{\bfseries \emph{Keywords:}}}


\section{Introduction}
\label{sec:introduction}

In this work, we study metric regularity properties of generalized equations of the form
\[ v \in F(x)\]
where $F : \Rn \rightrightarrows \Rn$ will correspond to the (limiting) subdifferential of a given function. Metric regularity of $F$ at $(\bar x,\bar v) \in \gph(F)$ is equivalent to the Aubin property of the local solution mapping $S : \Rn \rightrightarrows \Rn$, $S(v):=F^{-1}(v)$, at $(\bar v,\bar x)$: nearby solutions vary in a Lipschitzian set-valued manner under perturbations of $v$, cf.\ \cite[Chapter 9]{RockafellarWets}. Strong metric regularity requires that $S$ additionally admits a single-valued Lipschitz localization around $(\bar v,\bar x)$. Hence, the two notions can describe different branching properties: metric regularity permits several nearby solution branches, whereas strong metric regularity necessitates the existence of a unique branch.

For smooth square systems, the classical inverse function theorem makes metric and strong metric regularity equivalent. An analogous equivalence holds for important nonsmooth systems with polyhedral structure. In particular, the works by Robinson and Dontchev-Rockafellar show that, for smooth generalized equations over polyhedral convex sets, the Aubin property of the solution mapping forces strong regularity \cite{robinson1980regular,DontchevRockafellar}; see also \cite{robinson1992normal,ralph1993homeomorphism,ioffe2018polyhedral}. This polyhedral result has served as a long-standing benchmark for stability theory in variational analysis \cite{DontchevRockafellar,RockafellarWets}.

 Recent substantial progress has extended this ``polyhedral phenomenon'' to various structured non-polyhedral systems. In \cite{HangSarabi}, Hang and Sarabi verify that metric and strong metric regularity coincide under strict proto-differentiability and identify relative-interior conditions ensuring it for certain decomposable functions; see also \cite{OuyangMilzarekI,OuyangMilzarekII} for related results. In a parallel line of work, Chen, Chen, Sun, and colleagues investigate canonically perturbed KKT systems over non-polyhedral cones, \cite{chen2025aubin,chen2026characterization}. In \cite{chen2025aubin}, Chen et al.\ show the equivalence for nonlinear second-order cone programming at locally optimal solutions, while analogous results for nonlinear semidefinite programming are obtained in \cite{chen2026characterization}. Most recently, Ma and Sun \cite{MaSun} prove the equivalence for generalized equations over arbitrary closed convex $C^2$-cone-reducible sets, at general reference points and without the local-optimality assumption appearing in \cite{chen2025aubin,chen2026characterization}. Together, these results show that metric and strong metric regularity are equivalent for a broad class of non-polyhedral systems, but they also indicate that fundamental geometric structure is involved. 
 
 
This progress leads to the following natural question: Does this equivalence carry over to other rich classes of problems, such as, general \emph{definable} systems? In particular, let $h : \Rn \to \Rex$ be proper, lower semicontinuous, definable, prox-regular, and subdifferentially continuous at $\bar x$ for $\bar v \in \partial h(\bar x)$. Is tameness sufficient to guarantee
\begin{equation} \partial h \; \text{is metrically regular at $(\bar x,\bar v)$} \quad \implies \quad \partial h \; \text{is strongly metrically regular at $(\bar x,\bar v)$?} \label{eq:main-problem} \end{equation}
In general, the answer is negative. Ma and Sun construct a three-dimensional closed convex semialgebraic set $S$ and a smooth function $g$ such that $\nabla g + N_S = \partial(g + \iota_S)$ is metrically regular at a given reference pair $(\bar x,\bar v)$, while the associated localized inverse has two distinct solution branches \cite[Section 5]{MaSun}. The function $g+\iota_S$ is semialgebraic, prox-regular, and subdifferentially continuous at $(\bar x,\bar v)$. Thus, a definable or even semialgebraic structure does not, by itself, rule out branching. Moreover, the example by Ma and Sun can be adjusted to realize any prescribed finite number of local inverse branches, see \cite[Remark 5.4]{MaSun}. This motivates two complementary directions:
\begin{itemize}
\item Can we characterize and identify an additional local structure to ensure unique inverse branching at a general point $(\bar x,\bar v) \in \gph(\partial h)$? Strict proto-differentiability and $C^2$-cone reducibility provide two such mechanisms, but it is natural to seek conditions formulated intrinsically through the second-order geometry of $h$. 
\item Drusvyatskiy, Mordukhovich, and Nghia conjectured that if $\bar x$ is a local minimizer of a prox-regular and subdifferentially continuous function $h$, then metric regularity of $\partial h$ at $(\bar x,0)$ already implies strong metric regularity \cite[Conjecture 4.7]{DMN}. In view of the connections between strong metric regularity, tilt stability, and uniform quadratic growth \cite{DrusvyatskiyLewis}, this asks whether local minimality is sufficient to eliminate all competing stable branches.
\end{itemize}
The present paper addresses both directions in the definable setting. More precisely, we derive a point-wise, topological second-order criterion that characterizes the strong metric regularity of $\partial h$ at a metrically regular reference pair $(\bar x,\bar v)$. We then use this criterion to prove a strengthened version of the conjecture of Drusvyatskiy, Mordukhovich, and Nghia.

Let $h$ satisfy the above assumptions and let $(\bar x,\bar v)\in\gph(\partial h)$ be fixed. As a first intermediate result, we show that definability and prox-regularity imply that $h$ is properly twice epi-differentiable at $\bar x$ for $\bar v$. Under subdifferential continuity, its second subderivative $q:=\mathrm{d}^2h(\bar x|\bar v)$ then gives the exact graphical model
\begin{equation}
D(\partial h)(\bar x,\bar v)=\partial(\half q),
\label{eq:intro-proto-model}
\end{equation}
where $D(\partial h)$ denotes the proto-derivative of $\partial h$. The second subderivative $q$ is proper, lower semicontinuous, definable, positively two-homogeneous, and weakly convex, and its effective domain $\dom(q)$ is a convex cone. We further define 
\begin{equation}
\eq := \Span(\dom(q)), \quad \linq :=\dom(q)\cap(-\dom(q)), \quad \text{and} \quad \negq:=\{w\in\dom(q)\cap\Sph_{\eq}:q(w) <0\},
\label{eq:intro-negative-link}
\end{equation}
where $\Sph_{\eq} := \{x \in \eq: \|x\| = 1\}$. The set $\linq$ is the lineality space of $\dom(q)$ and the negative link $\negq$ contains all unit directions for which the second-order model $q = \mathrm{d}^2h(\bar x|\bar v)$ has negative curvature. Our first main result is an exact topological characterization of strong metric regularity based on $\negq$. \vspace{1mm}
\begin{mdframed}[
  linecolor=white!100,
  linewidth=.2mm,          
  backgroundcolor=gray!10, 
  roundcorner=0pt,         
  innerleftmargin=10pt,    
  innerrightmargin=15pt,   
  innertopmargin=7pt,     
  innerbottommargin=7pt   
]
\begin{theorem}[A topological curvature criterion for strong metric regularity, informal form]
\label{thm:intro-euler}
Assume (in addition) that $\partial h$ is metrically regular at
$(\bar x,\bar v)$.  Then, it holds that
\begin{equation}
\partial h\text{ is strongly metrically regular at }(\bar x,\bar v)
\quad\Longleftrightarrow\quad
\chi(\negq)\in\{0,2\}.
\label{eq:intro-euler-criterion}
\end{equation}
%
\end{theorem}
\end{mdframed}
Here, $\chi(\negq)$ denotes the Euler characteristic of the set $\negq$. The Euler characteristic is a number that describes the topological shape or structure of an underlying set; we refer to \cref{sec:intro-euler} for additional background and information. Hence,  \Cref{thm:intro-euler} expresses strong metric regularity as a topological property of the set of negative curvature directions of $h$. Although this criterion provides a complete characterization, it is topological and may not always be straightforward to verify. Our second main result is based on a simpler geometric sufficient condition involving only the directions in the lineality space $\linq$. \vspace{1mm}


\begin{mdframed}[
  linecolor=white!100,
  linewidth=.2mm,          
  backgroundcolor=gray!10, 
  roundcorner=0pt,         
  innerleftmargin=10pt,    
  innerrightmargin=15pt,   
  innertopmargin=7pt,     
  innerbottommargin=7pt   
]
\begin{theorem}[Lineality principle, informal form]
\label{thm:intro-lineality}
Assume that $\partial h$ is metrically regular at
$(\bar x,\bar v)$.  If
\begin{equation}
q(w) = \mathrm{d}^2h(\bar x|\bar v)(w) \geq0 \qquad \forall~w\in\linq,
\label{eq:intro-lineality-nonnegative}
\end{equation}
then there exists $\gamma >0$ such that  
\begin{equation}
 q(w) \geq \gamma \norm{w}^2 \qquad \forall~w\in\dom(q),
\label{eq:intro-uniform-positive}
\end{equation}
and $\partial h$ is strongly metrically regular at $(\bar x,\bar v)$.
\end{theorem}
\end{mdframed}

Thus, under metric regularity, nonnegative second-order curvature on the lineality space propagates to uniform quadratic positivity on the entire effective domain. In particular, the condition is automatic when $\dom(q)$ is pointed, i.e., $\linq = \{0\}$.
Property \eqref{eq:intro-lineality-nonnegative} is an intrinsic analogue of weak second-order necessary conditions used in constrained optimization. In nonlinear and conic programming, such conditions test the curvature of the Hessian of the Lagrangian, including an appropriate sigma term, on the lineality space of a critical cone; see, e.g., \cite[Definition~2.1]{AndreaniHaeserRamosSilva} or \cite[Definition~1]{FukudaHaeserMito}. Here, the effective domain $\dom(q)$ and its lineality space $\linq$ are determined directly by the second subderivative. In the following, we will call $\bar x$ a \emph{weak second-order stationary point} of the tilted function $h_{\bar v}:=h-\iprod{\bar v}{\cdot}$ when \eqref{eq:intro-lineality-nonnegative} holds. 

\Cref{thm:intro-euler} can be used to immediately establish the conjecture of Drusvyatskiy, Mordukhovich, and Nghia in the definable setting. If $\bar x$ is a (standard) \emph{second-order stationary point} of $h_{\bar v}$ in the sense that
\begin{equation}
q(w) = \mathrm{d}^2h_{\bar v}(\bar x|0)(w) \geq0 \qquad \forall~w\in\dom(q),
\label{eq:intro-second-order-stationarity}
\end{equation}
then we obtain $\negq = \varnothing$. Consequently, using $\chi(\varnothing)=0$,  \Cref{thm:intro-euler} yields
\begin{equation}
\left.
\begin{array}{c}
\bar x \; \text{is a second-order stationary point of} \; h_{\bar v},\\
\partial h \; \text{is metrically regular at} \; (\bar x,\bar v)
\end{array}
\right\}
\quad\implies\quad
\partial h \; \text{is strongly metrically regular at} \; (\bar x,\bar v).
\label{eq:intro-second-order-consequence}
\end{equation}
Every local minimizer of $h_{\bar v}$ satisfies \eqref{eq:intro-second-order-stationarity}, \cite[Theorem~13.24 (a)]{RockafellarWets}. Taking $\bar v = 0$ therefore proves the definable case of \cite[Conjecture 4.7]{DMN}. \Cref{thm:intro-lineality} strengthens this implication within the definable setting: local minimality, and even nonnegativity of $q$ on its full domain, can be replaced by the weaker lineality condition \eqref{eq:intro-lineality-nonnegative}.

The proof of \Cref{thm:intro-euler} is built on several steps. We first use a localized Moreau envelope-based transformation to replace the original subgradient inclusion by a $C^{1,1}$-gradient equation while preserving both metric and strong metric regularity. The semiderivative of the gradient of this localized Moreau envelope is the positively homogeneous gradient model
$H = \nabla\env{q/2}{\lambda}$,
where $\env{q/2}{\lambda}(z) := {\inf}_{y \in \Rn} \tfrac12 q(y) + \tfrac{1}{2\lambda}\|z-y\|^2$, $\lambda > 0$. We then construct a proximal homotopy to connect the Euler characteristics of the sets 
\[ \negq \qquad \text{and} \qquad \{z \in \Sph_{\eq}: \env{q/2}{\lambda}(z) < 0 \}, \] 
and
%
%
a homogeneous gradient-index formula converts the Euler characteristic into the local (Brouwer) degree $\degree_{\mathrm{loc}}(H,0)$ of $H$ at $0$. Specifically, we show that $\degree_{\mathrm{loc}}(H,0) = 1 - \chi(\negq)$. Metric regularity ensures that the definable mapping $H$ is open with finite fibers, and the degree-based condition $|\degree_{\mathrm{loc}}(H,0)| = 1$ then allows us to exclude multiple local inverse branches. The proof of the lineality principle is geometric rather than topological: We leverage (uniform) metric regularity to locally bound the negative inertia of the Hessian of the Moreau envelope (where it exists) across the smooth definable strata, while an induction on $\dim(\eq)$ is used to treat the case $\dom(q) \neq \eq$. The induction step provides a positive definite Hessian; inertia propagation then makes the Moreau envelope strongly convex. This will yield \eqref{eq:intro-uniform-positive} and strong metric regularity of $\partial h$.

In the following, we summarize our key contributions.

\begin{itemize}
\item \textbf{Second-order models}. We verify that definability and prox-regularity yield proper twice epi-differentiability; under subdifferential continuity, the second subderivative then gives the proto-derivative formula \eqref{eq:intro-proto-model}. Based on a localized Moreau envelope transform, we show that metric and strong metric regularity can be transferred between the original subgradient mapping $\partial h$ and several homogeneous models associated with $q = \mathrm{d}^2h(\bar x|\bar v)$.

\item \textbf{An Euler characteristic-based criterion}. We show that the topology of the set $\negq$ of negative curvature directions of $q$ provides a necessary-and-sufficient criterion for metric regularity to upgrade to strong metric regularity. This result combines various techniques and concepts from variational analysis (second subderivatives, open mappings), geometry (definability, stratifications), and topology (degree theory, Euler characteristic). 

\item \textbf{Lineality principle}. We verify that nonnegativity of $q$ on $\linq$ propagates, under metric regularity, to uniform quadratic positivity on $\dom(q)$. This proves the conjecture by Drusvyatskiy, Mordukhovich, and Nghia in the definable setting and strengthens its conclusion to weak second-order stationary points of tilted functions.
%
\end{itemize}

The remainder of the paper is organized as follows. \Cref{sec:prelim} presents several variational, geometric, and degree-theoretic preliminaries. \Cref{sec:model-reduction} establishes proper twice epi-differentiability and the exact Moreau reduction. \Cref{sec:euler} proves the Euler-link criterion and its second-order stationarity consequences. \Cref{sec:lineality} develops the lineality principle, and \Cref{sec:conclusion} concludes with open directions.

\section{Variational, definable, and topological preliminaries}
\label{sec:prelim}

We work in a fixed o-minimal expansion of the real field, and the word
\emph{definable} always refers to this structure.  Definable sets are closed
under finite Boolean operations, Cartesian products, and coordinate
projections, and every definable subset of $\R$ is a finite union of points
and intervals.  A mapping is definable when its graph is definable; in
particular, semialgebraic sets and mappings are definable.  We use
monotonicity for univariate functions \cite[Theorem~4.1]{dries1996geometric}, finite Whitney stratifications \cite[Theorem~4.8(1)]{dries1996geometric}, and definable triangulation \cite[Theorem~4.10]{dries1996geometric}. For further background on definability, we refer to \cite[Chapters~2--4]{vanDenDries}.

\subsection{Notation and variational tools}

\paragraph{General notation.} By $\iprod{\cdot}{\cdot}$ and $\|\cdot\|$,
we denote the Euclidean inner product and norm.  For matrices,
$\|\cdot\|$ is the spectral norm.  The sets of symmetric and symmetric
positive-definite $n\times n$ matrices are denoted by $\mathbb S^n$ and
$\mathbb S_{++}^n$, respectively.  For $A,B\in\mathbb S^n$, we write
$A\succeq B$ when $A-B$ is positive semidefinite.  The \emph{inertia} of
$A\in\mathbb S^n$ is
\[
\operatorname{In}(A):= (n_+(A),\nminus(A),n_0(A)),
\]
where the entries count the number of positive, negative, and zero eigenvalues of $A$,
respectively, \cite[Section~4.5]{HornJohnson2013}.  The
open ball centered at $x$ with radius $r$ is
$\bB_r(x):=\{y\in\R^n:\|y-x\|<r\}$; we write $\bB_r:=\bB_r(0)$. Furthermore, we use $\mathbb S_{n} := \{x \in \R^{n+1}: \|x\| = 1\}$ to denote the unit $n$-sphere in $\R^{n+1}$.

\paragraph{Sets.} We use the standard Painlev\'{e}-Kuratowski convergence for sets. For a sequence $\{C_k\}_k \subseteq \R^n$, the outer and inner limit are defined via 
\[ {\limsup}_{k\to\infty}~C_k
:=\{x\in\R^n:\exists\ k_j\uparrow\infty,\ \exists\ x_j\in C_{k_j},\
 x_j\to x\},
\]
and $\liminf_{k\to\infty} C_k:=\{x\in\R^n:~\exists k_0 \in \mathbb N, \forall k\geq k_0,~ \exists~x_{k}\in C_{k},~x_{k}\to x\}$. Furthermore, we write $C_k\to C$ when both limits coincide. 
The sets $\cl(C)$, $\co(C)$, and $C^\circ$ denote the closure, convex hull, and polar cone of $C \subseteq \Rn$. We use $\iota_C$ and $\sigma_C$ to denote the indicator and support function of $C$. For some fixed $x \in C$, the tangent and regular normal cone to $C$ at $x$ are defined by
\[ T_C(x) := {\limsup}_{t \downarrow 0}~t^{-1}({C-x}) \quad \text{and} \quad \hat N_C(x) := T_C(x)^\circ. \]
The (limiting/Mordukhovich) normal cone $N_C(x)$ to $C$ at $x$ is given by $N_C(x) := \limsup_{y\to x} \hat N_C(y)$. For $x \notin C$, we set $\hat N_C(x) = \varnothing$ and $N_C(x) = \varnothing$. When $C$ is convex, both normal cones coincide with the standard normal cone used in convex analysis.

\paragraph{Functions and derivatives.} The effective domain of a function $h : \Rn \to \Rexx$ is defined as $\dom(h):=\{x \in \Rn : h(x)<\infty\}$. We call $h$ proper if $h(x) > -\infty$ for all $x$ and $\dom(h) \neq \varnothing$. The epigraph $\epi(h)$ of a function $h:\R^n\to\Rex$ is given by $\epi(h):=\{(x,t)\in\R^n\times \R: h(x) \leq t\}$. A function $h:\R^n\to \Rexx$ is said to be lower semicontinuous (lsc) at $x$ if $\liminf_{y\to x} h(y)\geq h(x)$ and it is said to be lsc if it is lsc at every point. 
A sequence of extended real-valued functions $\{h_k\}_k$ epi-converges to $h$ if $\epi(h_k)\to\epi(h)$ (in the Painlev\'{e}-Kuratowski sense). We refer to \cite[Chapter 7]{RockafellarWets} for alternative characterizations and additional information.
For a function $h:\R^n\to \Rex$ and a point $x\in\dom(h)$, the {(limiting) subdifferential} is given by $\partial h(x) := \{v \in \Rn: (v,-1) \in N_{\epi(h)}(x,h(x)) \}$. We further set $\partial h(x)=\varnothing$ when $x\notin\dom(h)$.
The {second subderivative} $\mathrm{d}^2 h(x|v)(w)$ of $h$ at $x$ relative to $v$ in the direction $w \in \Rn$ is defined as
\[ \mathrm{d}^2h(x|v)(w) = \liminf_{t\downarrow 0, \, \tilde w \to w}~\Delta_t^2  \;\! h(x|v)(\tilde w), \quad \Delta_t^2  \;\! h(x|v)(w) := \frac{h(x+tw)-h(x)-t \cdot \iprod{v}{w}}{\half t^2}. \]
We say that $h$ is {twice epi-differentiable} at $x$ for $v \in \Rn$ if the second-order difference quotients $w \mapsto \Delta_{t_k}^2  \;\! h(x|v)(w)$ epi-converge to $\rd^2h(x|v)$ for every sequence $t_k \downarrow 0$. The function $h$ is said to be {properly twice epi-differentiable} if, in addition, $\rd^2h(x|v)$ is proper.

Let $S:\R^n\rightrightarrows \R^m$ be a set-valued mapping. The graph of $S$ is $\gph(S) := \{(x,y)\in \R^n\times \R^m: y\in S(x)\}$. We use $\dom(S) := \{x \in \Rn: S(x)\neq\varnothing\}$, to denote the domain of $S$. 
The graphical derivative of $S$ at $x$ for $v$, $(x, v) \in \gph(S)$, is the set-valued mapping $DS(x,v) : \Rn \rightrightarrows \R^m$ defined via $\gph(DS(x,v)) = T_{\gph(S)}(x,v)$. When the ``$\limsup$'' in the definition of $T_{\gph(S)}(x,v)$ is a full set limit, we say that $S$ is proto-differentiable at $x$ for $v$, cf. \cite[Section 8.H]{RockafellarWets}. A function $F: \Rn \to \R^m$ is called {directionally differentiable} at $x$ in the direction $w \in \Rn$ if the limit $\lim_{t \downarrow 0} \Delta_t \;\! F(x)(w)$ exists. Furthermore, $F$ is called {semidifferentiable} at $x$ in the direction $w$ if the limit $\lim_{t \downarrow 0, \tilde w \to w} \Delta_t \;\! F(x)(\tilde w)$ exists. We use $F^\prime(x;w)$ to denote both limits.  Whenever the semiderivative exists (for all directions), it is positively one-homogeneous as a function of its direction, i.e.,
\[
F^\prime(x;tw)=tF^\prime(x;w) \qquad \forall~t>0,\ w\in\Rn.
\]
In general, a mapping $h : \Rn \to \Rex$ is said to be positively $p$-homogeneous (or homogeneous of degree $p$), $p > 0$, if $h(tw) = t^p h(w)$ for all $t > 0$ and $w \in \Rn$.

\subsection{Metric regularity, openness, fibers, and branches} \label{sec:intro-metric}

\paragraph{Metric regularity and openness.} A set-valued mapping $S:\R^n\rightrightarrows\R^m$ is \emph{metrically regular} at $(\bar x,\bar v) \in\gph(S)$ if there exist neighborhoods $U$ of $\bar x$ and $V$ of $\bar v$, and a rate $\kappa>0$, such that
\begin{equation}\label{eq:def-mr}
\dist(x,S^{-1}(v)) \le \kappa\,\dist(v,S(x)) \qquad \forall~(x,v)\in U\times V.
\end{equation}
The mapping $S$ is \emph{strongly metrically regular} at $(\bar x,\bar v)$ if, in addition, $S^{-1}$ admits a single-valued Lipschitz localization around $(\bar v,\bar x)$.  Strong metric regularity always implies metric regularity.

The mapping $S:\R^n\rightrightarrows\R^m$ is said to have the Aubin property at $\bar x$ for $\bar v \in S(\bar x)$, if $\gph(S)$ is locally closed at $(\bar x,\bar v)$ and there are neighborhoods $U$ of $\bar x$ and $V$ of $\bar v$, and $\kappa > 0$ such that 
\[ S(x^\prime) \cap V \subset S(x) + \kappa \|x - x^\prime\|\,\bB_1 \quad \forall~x, x^\prime \in U. \]
If $\gph(S)$ is locally closed at $(\bar x,\bar v)$, then $S$ is metrically regular at $(\bar x,\bar v)$ if and only if $S^{-1}$ has the Aubin property at $\bar v$ for $\bar x$.

Metric regularity of $S$ at $(\bar x,\bar v)$ (at a rate $\kappa$) is further known to be equivalent to \emph{linear openness} of $S$ at $(\bar x,\bar v)$ (at a rate $\kappa^{-1}$). In particular, there exists a neighborhood $U \times V$ of $(\bar x,\bar v)$ and $r > 0$ such that for all $(x,v) \in \gph(S) \cap U \times V$ and $\veps \in [0,r]$, we have 
\[ \bB_{\veps\kappa^{-1}}(v) \subset S(\bB_{\veps}(x)), \]
cf. \cite[Theorem~3E.8]{dontchev2009implicit}. A function $F: \Rn \to \R^m$ is said to be \emph{open} (in a topological sense) if it maps open subsets of $\Rn$ to open subsets of $\R^m$. $F$ is \emph{locally open} on a neighborhood $U \subset \Rn$ if every sufficiently small neighborhood of $x \in U$ is mapped to a neighborhood of $F(x)$. Linear openness (of a single-valued mapping) clearly implies local openness. We refer to \cite{RockafellarWets,dontchev2009implicit} for more background.

\paragraph{Fibers.} Let $S:\R^n\rightrightarrows\R^m$ be a set-valued mapping. For $v\in \R^m$, the inverse image $S^{-1}(v) =\{x\in \Rn:v\in S(x)\}$ is also called the \emph{fiber} of $S$ over $v$. If $ U\subseteq \Rn$, then $S^{-1}(v)\cap U$ is a \emph{local fiber} over $v$. We say that the local fiber is \emph{finite} if the set $S^{-1}(v)\cap U$ is finite. We say that it is \emph{discrete} if each of its points is isolated in that fiber. Finite fibers are automatically discrete. 

\paragraph{Branches.} Let $ S : \Rn \rightrightarrows \R^m$ be set-valued and let $ U\subseteq \Rn$, $V\subseteq \R^m$ be given. A \emph{solution branch} of the generalized equation $v \in S(x)$ over $V$ relative to $U$ is a single-valued selection $s: V\to U$ satisfying $s(v) \in S^{-1}(v)$, $v \in V$. If $F : \Rn \to \R^m$ is single-valued, then a branch of $F^{-1}$ is called an \emph{inverse branch}. Strong metric regularity of $S$ at $(\bar x,\bar v)$ means that, after localization near $(\bar v,\bar x)$, the inverse $S^{-1}$ has exactly one solution branch and this branch is Lipschitz.

\subsection{Brouwer degree}

The Brouwer degree is a fundamental tool for counting local inverse branches with algebraic multiplicity, and will play a central role in our analysis. 
We now record several definitions and standard facts, following the
monograph \cite[Section~1.2]{DincaMawhin}.

\begin{definition}[Brouwer degree]\label{def:brouwer-degree}
The triple $(F,\Omega,y)$ is called \emph{admissible} if $\Omega\subset\Rn$ is bounded and open, $F : \cl(\Omega)\to\R^n$ is continuous, and $y\notin F(\partial\Omega)$. The Brouwer degree is the unique integer-valued map $(F,\Omega,y)\mapsto\degree(F,\Omega,y)$ on admissible triples satisfying the following axioms.
\begin{enumerate}[label=\textup{(D.\arabic*)},leftmargin=3.2em]
\item \label{D1} \emph{Normalization.}
If $y\in\Omega$, then $\degree(\mathrm{id},\Omega,y)=1$.

\item \label{D2} \emph{Additivity.}
If $\Omega_1,\Omega_2\subset\Omega$ are disjoint open sets and $y\notin F(\cl(\Omega)\setminus (\Omega_1\cup\Omega_2))$, then
\[
\degree(F,\Omega,y) =\degree(F,\Omega_1,y)+\degree(F,\Omega_2,y).
\]

\item \label{D3} \emph{Homotopy invariance.} If $\mathscr H : \cl(\Omega)\times[0,1]\to\R^n$ and $q : [0,1]\to\R^n$ are continuous and $q(t)\notin \mathscr H(\partial\Omega,t)$ for every $t\in[0,1]$, then $\degree(\mathscr H(\cdot,t),\Omega,q(t))$ is independent of $t$.
\end{enumerate}
\end{definition}

This is the axiomatic characterization of the Brouwer degree found in
\cite[Theorem~3.9 and Section~3.4]{BenevieriFuriPeraSpadini}; see also
\cite[Section~1.2]{DincaMawhin}.
The following \emph{excision property} is an immediate consequence of \ref{D2}.  If $V\subset\Omega$ is open and $F^{-1}(y)\cap\Omega\subset V$, then $(F,V,y)$ is admissible and
\begin{equation*}\tag{D.4}\label{eq:degree-excision}
\degree(F,\Omega,y)=\degree(F,V,y).
\end{equation*}
Indeed, applying \ref{D2} first with $\Omega=\Omega_1=\Omega_2=\varnothing$ gives $\degree(F,\varnothing,y)=0$, and applying it again with $\Omega_1=V$ and $\Omega_2=\varnothing$ gives \eqref{eq:degree-excision}; see~\cite[Theorem~3.10]{BenevieriFuriPeraSpadini} or \cite[Section 1.2.5]{DincaMawhin}.

\begin{definition}[Local Brouwer degree]\label[definition]{def:local-degree} Let $U\subset\R^n$ be open, let $F: U\to\R^n$ be continuous, and suppose that $x\in U$ is an isolated point of $F^{-1}(F(x))$. An \emph{isolating neighborhood} of $x$ is a bounded open set $V\subset\R^n$ such that $x\in V$, $\cl(V) \subset U$, $F^{-1}(F(x))\cap\cl(V) = \{x\}$. The \emph{local Brouwer degree} of $F$ at $x$ is defined as 
\[ \degree_{\mathrm{loc}}(F,x):=\degree(F,V,F(x)). \]
\end{definition}

The local Brouwer degree is independent of the isolating neighborhood. Moreover, if $F$ is differentiable at $x$ and $DF(x)$ is nonsingular, then
\begin{equation*}
\degree_{\mathrm{loc}}(F,x)
=\operatorname{sgn}(\det(DF(x))).
\end{equation*}
In the following, we state standard consequences of the axioms and classical computational results that will be used in our analysis. First, the condition $\degree(F,\Omega,y)\ne0$ implies $y\in F(\Omega)$. Second, if $F$ is continuous and $F^{-1}(y)\cap\Omega$ is finite, then additivity gives $\degree(F,\Omega,y) = {\sum}_{x\in F^{-1}(y)\cap\Omega}\,\degree_{\mathrm{loc}}(F,x)$. These are, respectively, Theorems~3.11 and~3.15 in \cite{BenevieriFuriPeraSpadini}.
Moreover, let $A, B \in \R^{n \times n}$ be invertible matrices, and let $y \notin AF(\partial\Omega)$ be given. If $F$ is continuous on $\cl(\Omega)$, then the mapping $A \circ F \circ B$ is continuous on $\cl(B^{-1}(\Omega))$ and, according to \cite[Lemma 1.3.1]{DincaMawhin}, it holds that
\begin{equation} \label{eq:degree-linear-iso} \degree(A \circ F \circ B, B^{-1}(\Omega),y) = \sgn(\det(AB)) \cdot \degree(F,\Omega,A^{-1}y). \end{equation}
Let $\sL$ be a finite-dimensional Euclidean space, let $F:\sL\to \sL$ be
continuous, and suppose that $F(0)=0$ and that $0$ is isolated in
$F^{-1}(0)$. For sufficiently small $r>0$, we will use the abbreviation
\begin{equation*}
\deg_0(F)
:=\deg_{\mathrm{loc}}(F,0)
=\deg(F,\bB_r\cap \sL,0).
\end{equation*}
(As just noted and as in \Cref{def:local-degree}, excision makes this integer independent of the radius $r$). 
When $\sL=\{0\}$, we use the standard zero-dimensional convention
$\deg_0(F)=1$.

\paragraph{The open- and homological inverse-mapping theorem.} In the following, we state the definable open-mapping theorem of Dinh and Pham \cite[Theorem~3.1]{DinhPham} and a definable version of the homological inverse-mapping theorem \cite[Theorem~1.2]{BarretoFenilleHartmann}; see also \cite[Lemma~2.4]{MaSun}.

\begin{theorem}[Definable openness]
\label{thm:definable-openness}
Let $F: \Omega \to\Rn$ be continuous and definable on the connected open set $\Omega \subset \Rn$.
If $F$ is open, then every fiber is finite and it holds that 
\[
\det(DF(x))\geq0\quad\text{for all }x\in D_F,
\qquad\text{or}\qquad
\det(DF(x))\leq0\quad\text{for all }x\in D_F,
\]
where $D_F := \{x \in \Omega: \text{$F$ is differentiable at $x$}\}$. 
\end{theorem}

\begin{theorem}[Inverse-mapping theorem]
\label{lem:open-degree-one}
Let $F:\Omega\to\Rn$ be continuous, definable, and open on a connected open
set $\Omega\subset\Rn$.  Then every $x\in\Omega$ is isolated in its fiber,
so $\deg_{\mathrm{loc}}(F,x)$ is well-defined.  If
$|\deg_{\mathrm{loc}}(F,x)|=1$, then $F$ is a local homeomorphism at $x$.
\end{theorem}

%
%

\subsection{Homotopy and Euler characteristic} \label{sec:intro-euler}

By definable triangulation \cite[Theorem~4.10]{dries1996geometric}, every definable set $S$ has finite Betti numbers $b_j(S)$ and the homotopy type of a finite CW complex. In particular, for such a set $S$, we will work with the ordinary, homotopy-invariant Euler characteristic, not the compactly supported Euler characteristic. Using singular homologies with rational
coefficients, the Betti numbers and Euler characteristic can be defined via
\begin{equation}
b_j(A):=\dim_{\mathbb Q}(H_j(A;\mathbb Q)),
\qquad \chi(A):= {\sum}_{j\ge0}(-1)^j b_j(A),
\qquad \chi(\varnothing):=0.
\label{eq:homotopy-euler-definition}
\end{equation}
The sum is finite. The homological formula for the Euler characteristic
\cite[Theorem~2.44]{Hatcher} identifies this number with the alternating cell
count of a finite CW model and shows that it depends only on the homotopy type. { Every nonempty contractible set has Euler characteristic $1$; examples include a point, a ball, and an open interval. Moreover, $\chi(\mathbb S_k)=1+(-1)^k$ and $\chi(\mathbb S_1\times\mathbb S_1)=0$. Thus, a circle has Euler characteristic $0$, a two-sphere has Euler characteristic $2$, and a torus has Euler characteristic $0$. The Euler characteristic is additive over finite disjoint unions of connected components and satisfies $\chi(A\times B)=\chi(A)\chi(B)$ \cite[Section~2.2, Exercises~20 and 21]{Hatcher}. It is not generally additive over arbitrary definable partitions under the convention used here.}

Two continuous maps $f_0,f_1:A\to B$ are \emph{homotopic} if there exists a
continuous map $\mathscr H:A\times[0,1]\to B$ satisfying
$\mathscr H(a,0)=f_0(a)$ and $\mathscr H(a,1)=f_1(a)$ for every $a\in A$; we
then write $f_0\simeq f_1$.  Two spaces $A$ and $B$ are \emph{homotopy
equivalent} if there are continuous maps $f:A\to B$ and $g:B\to A$ such that
\[
g\circ f\simeq\operatorname{id}_A
\qquad\text{and}\qquad
f\circ g\simeq\operatorname{id}_B.
\]
See \cite[Chapter~0, p.~3]{Hatcher}. {Homotopy equivalent definable sets have isomorphic singular homology groups \cite[Corollary~2.11]{Hatcher}, and therefore equal Betti numbers and Euler characteristics by \cite[Theorem~2.44]{Hatcher}.}

Let $\mathbb L$ be a Euclidean space, let $\Sph_{\!\;\mathbb L}:=\{z\in \mathbb L:\norm{z}=1\}$, and let $C\subseteq \mathbb L$ be a positive cone. The radial map
\[
C\setminus\{0\}\rightarrow
(C\cap\Sph_{\!\;\mathbb L})\times(0,\infty),
\qquad
x \mapsto ({x}/{\norm{x}},\norm{x}),
\]
is a homeomorphism with inverse $(\omega,r)\mapsto r\omega$. Homotopy invariance and the product formula for the Euler characteristic then give
\begin{equation}
\chi(C\setminus\{0\})
=\chi(C\cap\Sph_{\!\;\mathbb L})\chi((0,\infty))
=\chi(C\cap\Sph_{\!\;\mathbb L}).
\label{eq:cone-link-euler}
\end{equation}
For a positively $2$-homogeneous function $f : \sL \to \R$, the set $\{z: f(z)<0\}$ is a positive cone that does not contain the origin. Hence, $\{z: f(z)<0\}$ and the negative link $\{z \in \sphL: f(z)<0\}$ have the same Euler characteristic.

\section{A Moreau envelope-based reduction approach}
\label{sec:model-reduction}

In this section, we introduce a local transformation between
the original subdifferential $\partial h$ and a gradient mapping based on the Moreau envelope. We then provide connections between the semiderivative of this gradient mapping and the Moreau envelope of the second subderivative.

\subsection{Basic proximal properties} \label{sec:basicprox}
We first state several basic properties of the proximity operator, the Moreau envelope, and its localized versions. For a proper, lsc function $h:\Rn \to \Rex$ and $\lambda > 0$, the proximity operator and the Moreau envelope of $h$ are defined as
\[\proxh : \Rn \rightrightarrows \Rn, \quad \proxh(x)\,:=\,\argmin_{y \in \Rn}~h(y) + \frac{1}{2\lambda}\|y-x\|^2 \]
and $\envh(x):= \inf_{y \in \Rn} h(y) + \frac{1}{2\lambda}\|y-x\|^2$, see \cite{moreau1965proximite} and \cite[Definition~1.22]{RockafellarWets}. Next, we introduce several additional regularity concepts for the mapping $h$.

We say that $h$ is \emph{prox-bounded}, if there exist $x\in \Rn$ and $\rho>0$ such that $\inf_{y\in\R^n} h(y)+\frac{\rho}{2}\|y-x\|^2>-\infty$. Let $h:\Rn \to \Rex$ be proper and lsc, and let $(\bar x,\bar v)\in\gph(\partial h)$ be given. The function $h$ is said to be \emph{prox-regular at $\bar x$ for $\bar v$} if there are $\veps>0$ and $\rho\ge0$ such that
\begin{equation}\label{eq:def-prox-regular}
h(x')\ge h(x)+\langle v,x'-x\rangle-\frac{\rho}{2}\|x'-x\|^2
\end{equation}
whenever $\|x'-\bar x\|<\veps$, $\|x-\bar x\|<\veps$, $\|v-\bar v\|<\veps$, $v\in\partial h(x)$, and $h(x)<h(\bar x)+\veps$. Furthermore, $h$ is \emph{subdifferentially continuous at $\bar x$ for $\bar v$} if
\[
(x^k,v^k)\to(\bar x,\bar v),\quad v^k\in\partial h(x^k) \quad \implies\quad h(x^k)\to h(\bar x).
\]
These are standard definitions; cf.\ \cite[Definitions~1.23, 13.27, and 13.28]{RockafellarWets}. The following well-known result is taken from \cite[Proposition 13.37]{RockafellarWets}; see also \cite[Theorem 4.4]{PolRoc96}.

\begin{proposition} \label[proposition]{prop:prox-env} Let $h : \Rn \to \Rex$ be lsc, prox-bounded, and prox-regular
at $\bar x$ for $\bar v\in\partial h(\bar x)$. Then there exist $\veps, \rho > 0$ such that for all $\lambda \in (0,1/\rho)$, there is a neighborhood $Z_\lambda$ of  $\bar z := \bar x+\lambda \bar v$ on which the following statements hold:
\begin{enumerate}[label=\textup{\textrm{(\roman*)}},topsep=2pt,itemsep=0ex,partopsep=0ex]
\item The proximity operator $\proxh$ is nonempty, single-valued, and Lipschitz continuous (on $Z_\lambda$) and we have $\proxh(\bar z) = \bar x$. 
\item The Moreau envelope $\envh$ is of class $C^{1,1}$ on $Z_\lambda$ and satisfies
$ \nabla\envh(z) = \frac{1}{\lambda}(z-\proxh(z)) = (\lambda I + T_{\veps}^{-1})^{-1}(z)$ ,
where $T_\veps$ is an $h$-attentive $\veps$-localization of $\partial h$ around $(\bar x,\bar v)$, i.e., $T_\veps(x) := \partial h(x) \cap \bB_{\veps}(\bar v)$ if $x \in \bB_{\veps}(\bar x)$ and $h(x) < h(\bar x) + \veps$ and $T_\veps(x) := \varnothing$ otherwise.
\end{enumerate}
\end{proposition}

If $h$ is subdifferentially continuous at $\bar x$ for $\bar v$, then the $h$-attentive localization $T_\veps$ in \Cref{prop:prox-env}~(ii) can be chosen as an ordinary localization, i.e., we can drop the condition ``$h(x) < h(\bar x) + \veps$''. We now introduce a localized Moreau envelope-based transform that allows us to transfer metric regularity properties of $\partial h$ to the somewhat simpler gradient of the Moreau envelope. Such localized constructions and techniques have appeared frequently in variational analysis; see \cite{PolRoc96,RockafellarWets}.
 
\begin{lemma}[Local Moreau transform and definability] \label[lemma]{lem:moreau-transform}
Let $h:\Rn\to\Rex$ be proper and lsc and let $r > 0$ be given. Suppose that $h$ is prox-regular 
at $\bar x$ for $\bar v \in \partial h(\bar x)$ and define the localized function 
\[ \tilde h:=h+\iota_{\cl(\bB_r(\bar x))}. \] 
There then exists $\lambda,\veps > 0$, a neighborhood $Z$ of $\bar z := \bar x + \lambda \bar v$, and an $h$-attentive $\veps$-localization $T_\veps$ of $\partial h$ (around $(\bar x,\bar v)$) such that the following properties are satisfied: 
\begin{enumerate}[label=\textup{\textrm{(\roman*)}},topsep=2pt,itemsep=0ex,partopsep=0ex]
\item The Moreau envelope $\envth$ is $C^{1,1}$ on $Z$ with $\nabla \envth(z) = (\lambda I + T_\veps^{-1})^{-1}(z)$  and $\proxth$ is nonempty, single-valued, and Lipschitz continuous on $Z$. 

\item The localization $T_\veps$ is (strongly) metrically regular at $(\bar x,\bar v)$ if and only if $\nabla \envth$ is (strongly) metrically regular at $(\bar z,\bar v)$.
\item If, in addition, $h$ is definable, then $\envth$, $\proxth$, and $\nabla \envth$ are definable and $\proxth$ and $\nabla \envth$ are semidifferentiable at $\bar z$ (in all directions). 
\end{enumerate}
\end{lemma}

\begin{proof} The function $\tilde h$ is proper, lsc, and its effective domain lies in the compact set $B := \cl(\bB_r(\bar x))$.
Thus, $\tilde h$ is bounded from below and hence prox-bounded.
Since $\bar x$ lies in the interior of $B$, $\tilde h$ agrees with $h$ on a neighborhood $U$ of the point $\bar x$.  In particular,
$\bar v \in\partial\tilde h(\bar x)$, $\tilde h$ is prox-regular
at $\bar x$ for $\bar v$, and $\partial \tilde h(x) = \partial h(x)$ for all $x \in U$.
Consequently, part~\textup{(i)} follows directly from \Cref{prop:prox-env}.
Moreover, after potentially shrinking $U$ and by setting $V:=\mathbb B_\epsilon(\bar v)$, we obtain
\[
\nabla \envth(x+\lambda v) = v \;\; \iff \;\; x + \lambda v \in (\lambda I + T_\veps^{-1})(v) \;\; \iff \;\; x\in T_\veps^{-1}(v) \;\; \iff \;\; v\in T_\veps(x)
\]
for all $(x,v)$ in the neighborhood $U \times V$ of $(\bar x,\bar v)$. Set $\Psi(v,x):=(v,x+\lambda v)$. Locally around
$(\bar v,\bar z)$, we then have
\[ \gph((\nabla \envth)^{-1}) =\Psi(\gph(T_\varepsilon^{-1})). \]
Thus, the inverse mappings differ locally by the Lipschitz translation $v\mapsto\lambda v$, which preserves the Aubin property and single-valued Lipschitz localizations.
Noticing that $\gph(T_\veps)$ is locally closed around $(\bar x,\bar v)$ (provided that $\veps > 0$ is sufficiently small; cf.\ \cite[Lemma 3.4]{gfrerer2025regular}), the inverse characterizations of metric and strong metric regularity therefore prove \textup{(ii)}. 

If $h$ is definable, then so is $\tilde h$, and the graphs of $\envth$ and $\proxth$ are definable. More specifically, $\envth$ is the marginal function of the definable mapping $e(x,y) := \tilde h(y) + \frac{1}{2\lambda} \|y-x\|^2$ and we can write $\gph(\proxth) = \{(x,y) : e(x,y) = \envth(x)\}$. Using the extensive calculus available for tame mappings, this implies definability of $\envth$ and $\proxth$ and hence, $\nabla \envth$ is definable as well. We refer to \cite[Chapter 1]{vanDenDries}, \cite[Section 1.3]{Coste}, and \cite[Section 3.3]{milzarek2016numerical} for further details.
As noted in \cite{bolte2009semismooth}, every definable locally Lipschitz mapping admits directional derivatives. (This follows from the monotonicity lemma \cite[Theorem 4.1]{dries1996geometric} applied to the definable first-order difference quotient.) Moreover, local Lipschitz continuity upgrades this to semidifferentiability. As a result, $\proxth$ and $\nabla \envth$ are semidifferentiable at $\bar z$.
\end{proof}

If $h$ is prox-bounded, then all results in \Cref{lem:moreau-transform}~(iii) hold for the original function $h$. Next, we show that definability --- in addition to the properties listed in \Cref{lem:moreau-transform}~(iii) --- also ensures proper twice epi-differentiability.

\begin{proposition}[Proper twice epi-differentiability via definability]
\label[proposition]{prop:automatic-twice-epi}
Let $h:\Rn\to\Rex$ be proper, lsc, and definable and suppose that $h$ is prox-regular at $\bar x$ for $\bar v \in \partial h(\bar x)$. Then, $h$ is properly twice epi-differentiable at $\bar x$ for $\bar v$. Furthermore, if  $h$ is subdifferentially continuous at $\bar x$ for $\bar v$, then $\partial h$ is proto-differentiable at $\bar x$ for $\bar v$ and we have
\begin{equation}
D(\partial h)(\bar x,\bar v)
=\partial(\half \rd^2h(\bar x|\bar v)).
\label{eq:automatic-proto}
\end{equation}
\end{proposition}
\begin{proof} Let $\tilde h$, $\lambda, \veps, r > 0$, and $T_\veps$ be given as in \Cref{lem:moreau-transform}. Since $\proxth$ is semidifferentiable at $\bar z$ (by \Cref{lem:moreau-transform}~(iii)), we can apply \cite[Exercises 13.18, 13.35, 13.45 and Theorem 13.40]{RockafellarWets} in tandem to infer that $\tilde h$ is twice epi-differentiable at $\bar x$ for $\bar v$ and $T_\veps$ is proto-differentiable at $\bar x$ for $\bar v$ with
\[  DT_\veps(\bar x,\bar v) = \partial(\tfrac12 \rd^2 \tilde h(\bar x|\bar v)). \]
Since $\tilde h$ agrees with $h$ in a neighborhood of $\bar x$, we have $\Delta_t^2  \;\! \tilde h(\bar x|\bar v) = \Delta_t^2  \;\! h(\bar x| \bar v)$ on each bounded set of directions and for all sufficiently small $t > 0$. This readily implies that the (second-order) epi-limits of $\tilde h$ and $h$ coincide and twice epi-differentiability is preserved. Moreover, the epi-limit $\rd^2 h(\bar x|\bar v)$ is proper. In particular, for some $\rho > 0$, the prox-regularity of $h$ (at $\bar x$ for $\bar v$) implies
\begin{equation} \label{eq:d2-prox-reg}
\Delta_t^2  \;\! h(\bar x|\bar v)(\tilde w)\geq -\rho \|\tilde w \|^2,
\end{equation}
on bounded sets of directions and for small $t$. Thus, it holds that $\rd^2 h(\bar x|\bar v)(w) \neq -\infty$ for all $w$.  At the origin, the constant test direction $\tilde w = 0$ gives $\rd^2 h(\bar x|\bar v)(0) \leq 0$ and the estimate \eqref{eq:d2-prox-reg} gives the reverse inequality. Hence, it follows $\rd^2h(\bar x|\bar v)(0)=0$.

Under subdifferential continuity, we can directly apply \cite[Theorem~13.40]{RockafellarWets} to obtain proto-differentiability of $\partial h$ and \eqref{eq:automatic-proto}.
\end{proof}

The second subderivative $q := \rd^2 h(\bar x|\bar v)$ is lsc and positively homogeneous of degree $2$. Moreover, if $h$ is prox-regular and properly twice epi-differentiable at $\bar x$ for $\bar v$, then for some $\rho\geq0$, the mapping 
\begin{equation}
w \mapsto q(w) + \rho \|w\|^2 = \rd^2 h(\bar x|\bar v)(w) + \rho \|w\|^2 \quad \text{is convex};
\label{eq:q-weak-convex}
\end{equation}
see \cite[Propositions~13.5 and~13.49]{RockafellarWets}, i.e., $q$ is $2\rho$-weakly convex. In particular, $\dom(q) = \dom(\rd^2 h(\bar x|\bar v))$ is a convex cone. (This cone need not be closed.) Moreover, if $h$ is definable, then the second subderivative $q$ is definable as well. (The graph of $q$ can be written as a ``first-order formula'' only involving definable operations \cite{Coste}.)

\subsection{Regularity via subderivative- and Moreau envelope-based models} \label{sec:regularity-moreau}

In this subsection, we further connect the localized result in \Cref{lem:moreau-transform}~(ii) to regularity properties of the two models
\[ H := \nabla \env{q/2}{\lambda} \quad \text{and} \quad M := \partial (\half q) \quad \text{where} \quad q := \rd^2 h(\bar x|\bar v), \quad \lambda > 0. \]
Let $h$ be properly twice epi-differentiable at $\bar x$ for
$\bar v \in \partial h(\bar x)$ and let $\tilde h := h + \iota_{\cl(\bB_r(\bar x))}$, $r > 0$, be the localization defined in \Cref{lem:moreau-transform}. Following the proof of \Cref{prop:automatic-twice-epi}, $\tilde h$ is then also twice epi-differentiable (at $\bar x$ for $\bar v$) with $\rd^2 \tilde h(\bar x|\bar v) = \rd^2 h(\bar x|\bar v)$.  Fix a sufficiently small $\lambda>0$ and set $\bar z:=\bar x+\lambda\bar v$.  Combining \cite[Exercises 13.18, 13.35, and 13.45]{RockafellarWets}, 
we can infer, 
for every $w\in\Rn$,
\begin{equation}
(\nabla\envth)^\prime(\bar z;w)
=\lambda^{-1} (w-\proxth^\prime(\bar z;w)) 
=\lambda^{-1} (w-\prox_{\lambda q/2}(w))
=\nabla \env{q/2}{\lambda}(w).
\label{eq:semi-grad-env-model}
\end{equation}
Since the second subderivative $q$ is a weakly convex mapping, the properties stated in \Cref{prop:prox-env} hold globally on $\Rn$ for $\env{q/2}{\lambda}$, $H = \nabla \env{q/2}{\lambda}$, and $\prox_{\lambda q/2}$ provided that $\lambda > 0$ is sufficiently small; see, e.g., \cite[Corollary~3.4]{HohLabObe20}. Furthermore, (as, e.g., shown in \cref{lem:moreau-transform} (iii)), definability of $h$ (and $q$) carries over to the models $H$ and $M$. In particular, $H$ is a definable, Lipschitz continuous, and positively homogeneous function.
  
Before presenting the main result of this section, we show that (strong) metric regularity generally passes to the semiderivative mapping.

\begin{lemma}[Regularity and semiderivatives]
\label[lemma]{lem:mr-passes-tangent}
Let $F:\R^n\to\R^n$ be locally Lipschitz and
semidifferentiable at $x$ with $F(x) = v$. If $F$ is (strongly) metrically regular at $(x,v)$, then the semiderivative mapping $w \mapsto F^\prime(x;w)$ is (strongly) metrically regular at $(0,0)$.
\end{lemma}

\begin{proof}
Suppose $F$ is metrically regular at $({x}, {v})$ with rate $\kappa$ and set $G := F^\prime(x;\cdot)$. Let $u, w \in \Rn$ be arbitrary. We aim to show $\dist(u, G^{-1}(w)) \leq \kappa \|G(u) - w\|$. For $t > 0$ (sufficiently small), let us define $x_t := {x} + tu$ and $v_t := {v} + tw$. By metric regularity, there exists $z_t \in F^{-1}(v_t)$ such that
\[ \|z_t - x_t\| \leq \kappa \|F(x_t) - v_t\| + o(t). \]
Setting $\zeta_t = (z_t - {x})/t$, we have $F({x} + t \zeta_t) = v_t$ and dividing the previous inequality by $t$, it follows:
\begin{equation} \label{eq:inter-01} \|\zeta_t - u\| \leq \kappa \left\| \frac{F({x} + tu) - F({x})}{t} - w \right\| + o(1). \end{equation}
As $t \downarrow 0$, semidifferentiability implies that the right hand-side of \eqref{eq:inter-01} converges to $\kappa \|G(u) - w\|$. 
Since $F$ is locally Lipschitz { and using \eqref{eq:inter-01}}, $\{\zeta_t\}_t$ is bounded. Let $\zeta^*$ be an accumulation point of $\{\zeta_t\}_t$; then we have $\|\zeta^* - u\| \leq \kappa \|G(u) - w\|$. Moreover, dividing $F({x} + t \zeta_t) -v = tw$ by $t$ and letting $t \downarrow 0$, we can infer $G(\zeta^*) = w$. Thus $\zeta^* \in G^{-1}(w)$, proving the first claim.

Suppose now that $F$ is strongly metrically regular at $x$ for $v = F(x)$. By the first part of the proof, $G$ is metrically regular (at $(0,0)$). In particular, $G^{-1}(w)\neq\varnothing$ for every $w\in\Rn$; hence, $G$ is surjective. Moreover, by \cite[Theorem 1F.2]{dontchev2009implicit}, there exist a neighborhood $U$ of $x$ and $\sigma > 0$ such that
\begin{equation} \label{eq:f-str-met} \sigma \|y - z\| \leq \|F(y) - F(z)\| \quad \forall~y, z \in U. \end{equation} 
Let $w_1, w_2 \in \R^n$ be fixed and set $y_t = x + tw_1$, $z_t = x + tw_2$ in \eqref{eq:f-str-met} for some $t > 0$. Dividing by $t$ and letting $t \downarrow 0$, the semidifferentiability of $F$ at $x$ then yields $\sigma \|w_1 - w_2\| \leq \|G(w_1) - G(w_2)\|$. 
Thus, $G$ is injective. This shows that $G$ is a bijection with a Lipschitz inverse. Hence, $w \mapsto F^\prime(x;w)$ is strongly metrically regular (at $(0,0)$).
\end{proof}

\begin{theorem}[Model-based reduction and regularity]
\label{thm:exact-reduction}
Let $h : \Rn \to \Rex$ be proper, lsc, and definable, and let
$(\bar x,\bar v)\in\gph(\partial h)$.  Assume that $h$ is prox-regular and
subdifferentially continuous at $\bar x$ for $\bar v$, and consider the
localized function $\tilde h:=h+\iota_{\cl(\bB_r(\bar x))}$, $r > 0$. If
$\partial h$ is metrically regular at $(\bar x,\bar v)$, then, for sufficiently small $\lambda > 0$ and setting $q := \rd^2 h(\bar x|\bar v)$, the following properties are equivalent:
\begin{enumerate}[label=\textup{(\roman*)}]
\item $\partial h$ is strongly metrically regular at $(\bar x,\bar v)$;
\item $\nabla \envth$ is strongly metrically regular at $(\bar z,\bar v)$ where $\bar z = \bar x + \lambda \bar v$;
\item $H=\nabla\env{q/2}{\lambda}$ is strongly metrically regular at
      $(0,0)$;
\item $M=\partial(\half q)$ is strongly metrically regular at $(0,0)$.
\end{enumerate}
All four model mappings in \textup{(i)}--\textup{(iv)} are metrically regular at their respective reference pairs.
\end{theorem}

\begin{proof}
By \Cref{prop:automatic-twice-epi}, the function $h$ is properly twice
epi-differentiable at $\bar x$ for $\bar v$, and its subdifferential $\partial h$ is proto-differentiable there. Since $h$ is assumed to be subdifferentially continuous at $\bar x$ for $\bar v$, the $h$-attentive $\veps$-localization $T_\veps$ introduced in \Cref{lem:moreau-transform} can be replaced by the subdifferential $\partial h$, i.e., we have $T_\veps^{-1}(v) \cap U = (\partial h)^{-1}(v) \cap U$ locally in a neighborhood $U \times V$ of $(\bar x, \bar v)$. Thus, the equivalence of (i) and (ii) and the metric regularity of $G$ follow from \Cref{lem:moreau-transform}~(ii). 

Using \eqref{eq:semi-grad-env-model} and the relation
\[ v \in \partial (q/2)(x) \quad \iff \quad \prox_{\lambda q/2}(x+\lambda v) = x, \]
cf.\ \cite[Proposition 3.1]{HohLabObe20}, we have $H(w + \lambda u) = u$ if and only if $u \in \partial (\frac12 q)(w) = M(w)$. This implies $H^{-1}(u) = M^{-1}(u) + \lambda u$. Since weakly convex functions are subdifferentially continuous, $\gph(M)$ is locally closed. We can now proceed as in the proof of \Cref{lem:moreau-transform}~(ii) to show that (strong) metric regularity of the functions $H$ and $M$ at $(0,0)$ is equivalent. Furthermore, the metric regularity of $H$ and $M$ and the direction ``(ii) $\implies$ (iii)'' follow from \Cref{lem:mr-passes-tangent}.

To complete the proof, we need to verify that (iii) implies (ii). Thus, suppose that $H$ is strongly metrically regular at $(0,0)$. Since $H$ is strongly metrically regular and positively homogeneous, its local inverse
globalizes to a positively homogeneous Lipschitz inverse on $\R^n$. Thus, $H$ is a global homeomorphism and we can infer $|\deg_0(H)|=1$ where $\deg_0 (H) := \deg_{\mathrm{loc}}(H,0)$. Moreover, it follows
\[ \eta := {\min}_{\|{w}\|=1}~\|{H(w)} \|>0 .
\]
Set $G:=\nabla \envth$ and define $\mathscr H : \cl(\bB_1) \times \R_+ \to \R^n$ where $\mathscr H(w,t) := \Delta_t \;\!G(\bar z)(w)$ if $t > 0$ and $\mathscr H(w,0) := G^\prime(\bar z;w) = H(w)$. Since the homotopy map $w \mapsto \mathscr H(w,t)$ converges locally uniformly to $H$ as $t\downarrow0$ (by the local Lipschitz continuity and semidifferentiability of $G$), there exists $\delta > 0$ such that
\begin{equation}
{\sup}_{\|{w}\|=1}~\|{\mathscr H(w,t)-H(w)} \| <\frac{\eta}{2} \quad \forall~t \in [0,\delta]
\label{eq:uniform-degree-bound}
\end{equation}
and $\mathscr H$ is continuous on $\cl(\bB_1) \times [0,\delta]$. This implies that $0 \notin \mathscr H(\partial\bB_1,t)$ for all $t \in [0,\delta]$. Moreover, after shrinking to a connected definable neighborhood $\Omega$ of $\bar z$, metric regularity makes $G|_\Omega$ open (at a
linear rate) at every point. The definable open-mapping theorem
\cite[Theorem~3.1]{DinhPham} (see also \Cref{lem:open-degree-one}) then gives finite fibers, so $\bar z$ is isolated in its fiber $G^{-1}(\bar v)$. Thus, applying the homotopy invariance of the Brouwer degree, we can infer
\begin{equation}
\deg_{\mathrm{loc}}(G,\bar z) \underset{\text{(a)}}{=} \deg(\mathscr H(\cdot,t),\bB_1,0) = \deg_0(H) = \pm 1, \quad t \in (0,\delta].
\label{eq:degree-G-H}
\end{equation}
More specifically, to obtain (a), we may set $\mathscr H(w,t) = t^{-1} g(tw)$ where $g(w) := G(\bar z + w)- G(\bar z)$. Using \eqref{eq:degree-linear-iso}, it then follows $\deg(\mathscr H(\cdot,t),\bB_1,0) = \deg((t^{-1}I)\circ g \circ (tI),\bB_1,0) = \deg(g,\bB_t,0) = \deg_{\mathrm{loc}}(G,\bar z)$. Consequently, by \Cref{lem:open-degree-one}, $G = \nabla \envth$ is a local homeomorphism. Metric regularity then ensures the (local) Lipschitz continuity of the inverse mapping $G^{-1}$ via the Aubin property. Thus (ii) holds.
\end{proof}

\Cref{thm:exact-reduction} allows existing regularity results to be applied to the second-order model $M = \partial(\frac12q)$ rather than directly to $h$ or $\partial h$. We briefly illustrate this observation in the following example.  

\begin{example}[Structured second-order models] \label[example]{ex:structured-models}
Suppose that the second subderivative $q$ has the generalized quadratic form
\begin{equation} \label{eq:example-quadratic}
    q(w)=\langle w,Qw\rangle+\iota_K(w),
\end{equation}
where $Q\in\mathbb S^n$ and $K$ is a linear subspace. Then $M(w)=Qw+N_K(w)$. Let $P$ have orthonormal columns spanning $K$, and set $A:=P^\top QP$. For $w=Py$, we have
\[
    v\in M(Py)\quad\Longleftrightarrow\quad
    P^\top v=Ay.
\]
Metric regularity forces $A$ to be nonsingular: otherwise $\operatorname{range}(M) =P\operatorname{range}(A)+K^\perp$ is a lower-dimensional linear subspace, contradicting local openness. Consequently, $M^{-1}(v)=\{PA^{-1}P^\top v\}$ is single-valued and Lipschitz. Hence $\partial h$ is strongly metrically regular at $(\bar x,\bar v)$. If $K=\{0\}$, the same conclusion follows from $M^{-1}(v)=\{0\}$.

Generalized quadratic models of the form \eqref{eq:example-quadratic} play an important role in variational analysis \cite{PolRoc96,poliquin1996hessian,hang2024strict,HangSarabi,OuyangMilzarekI} and are often linked to certain inherent smoothness properties. In particular, the presented example covers the strictly proto-differentiable setting, in which the second subderivative is generalized quadratic under the standing assumptions; see, e.g., \cite[Proposition~3.2]{HangSarabi}.
More generally, suppose that $K$ is a nonempty closed convex cone. If $K$ is polyhedral, the Dontchev--Rockafellar theorem \cite{DontchevRockafellar} implies that metric regularity of $M$ yields strong metric regularity.
The same conclusion follows from \cite[Theorem~4.4]{MaSun} when $K$ is $C^2$-cone reducible at each of its points. 
\Cref{thm:exact-reduction} therefore yields strong metric regularity of $\partial h$ at $(\bar x,\bar v)$ in such cases. In \Cref{sec:euler}, we explore an alternative topological criterion that characterizes strong metric regularity without additional structural conditions on the model $q$.
\end{example}

Finally, to conclude this section, we show that the different regularity properties discussed in \Cref{thm:exact-reduction} can be preserved when the mappings $H = \nabla\env{q/2}{\lambda}$ and $M = \partial(\frac12 q)$ are restricted to the linear space $\eq := \Span(\dom(q))$. 

Set $\dim(\eq) = r$ and let $P \in \R^{n \times r}$ and $Q \in \R^{n \times (n-r)}$ be orthonormal bases of $\eq$ and $\eq^\perp$. Let $\mathcal P_\eq := PP^\top$ and $\mathcal P_\perp := QQ^\top$ denote the orthogonal projectors onto $\eq$ and $\eq^\perp$, respectively. We define the restricted mappings
\begin{subequations}
\begin{align}
\label{eq:ortho-1} & q_\eq : \R^r \to \Rex, \quad q_\eq(y) := q(Py) = \mathrm{d}^2 h(\bar x|\bar v)(Py), \quad M_\eq : \R^r \rightrightarrows \R^r, \quad M_\eq(y) := \partial(\tfrac12 q_\eq)(y), \\
\label{eq:ortho-2} & e_\eq : \R^r \to \R, \quad e_\eq(y) := \env{q_\eq/2}{\lambda}(y), \quad \text{and} \quad H_\eq : \R^r \to \R^r, \quad H_\eq(y) := \nabla e_\eq(y). 
\end{align}
\end{subequations}
Setting $z_\eq := P^\top z \in \R^r$, $z_\perp := Q^\top z \in \R^{n-r}$, and $T = [P \;\, Q ]$, the envelope $\env{q/2}{\lambda}$ can be decomposed as follows
\begin{align*} \env{q/2}{\lambda}(z) & = \min_{w \in \eq}~\frac12 q(w) + \frac{1}{2\lambda} \|w-[\mathcal P_\eq z+\mathcal P_\perp z]\|^2 = \min_{w \in \R^n}~\frac12 q(\mathcal P_\eq w) + \frac{1}{2\lambda}\|P^\top [w - z]\|^2 + \frac{1}{2\lambda} \|Q^\top z\|^2 \\ & = \min_{y \in \R^r}~\frac12 q(Py) + \frac{1}{2\lambda} \|y-z_\eq\|^2 + \frac{1}{2\lambda}\|z_\perp\|^2 = e_\eq(z_\eq) + \frac{1}{2\lambda}\|z_\perp\|^2 \end{align*}
and we obtain
\begin{equation} \label{eq:ortho-gradient-splitting} T^\top H(z) = T^\top \nabla \env{q/2}{\lambda}(z) = T^\top P\nabla e_\eq(P^\top z) + \lambda^{-1} T^\top \mathcal P_\perp z = \begin{bmatrix} H_\eq(z_\eq) \\ \lambda^{-1}z_\perp \end{bmatrix}. \end{equation}
Combining \eqref{eq:degree-linear-iso}, the Cartesian-product formula for the Brouwer degree \cite[Lemma~7.1.1 and Theorem~7.1.1]{DincaMawhin}, and \ref{D1}, we can therefore infer
\begin{equation} \label{eq:deg-ortho} \deg_0(H) = \deg_0(T^\top \circ H \circ T)  =\deg_0(H_\eq)\cdot \deg_0(\lambda^{-1} \Id )
=\deg_0(H_\eq),
\end{equation}
since $\lambda > 0$ (provided that $0$ is isolated in $H^{-1}(0)$ or, equivalently, in $H_\eq^{-1}(0)$). Moreover, using \eqref{eq:ortho-gradient-splitting}, we can establish the following inverse-image identity
\begin{equation}
H^{-1}(u) = PH_\eq^{-1}(u_\eq) + \lambda Qu_\perp, \quad u_\eq = P^\top u,\quad u_\perp = Q^\top u.
\label{eq:ortho-inverse-splitting}
\end{equation}
Consequently, for $z = Pz_\eq+Qz_\perp$ and $u=Pu_\eq+Qu_\perp$, we obtain
\begin{equation*} 
\begin{split}
\dist(z,H^{-1}(u))^2
&=\dist(z_\eq,H_\eq^{-1}(u_\eq))^2
  +\norm{z_\perp-\lambda u_\perp}^2,\\
\norm{H(z)-u}^2
&=\norm{H_\eq(z_\eq)-u_\eq}^2
  +\lambda^{-2}\norm{z_\perp-\lambda u_\perp}^2.
\end{split}
\end{equation*}
These identities show directly that $H$ is metrically regular if and only if
$H_\eq$ is metrically regular. The inverse formula \eqref{eq:ortho-inverse-splitting} further shows that the local single-valuedness and Lipschitz continuity of the inverse localizations of $H$ and $H_\eq$ are equivalent; thus, strong metric regularity of $H$ and $H_\eq$ are also equivalent. 

Similar results and connections can be established for $q_\eq$ and $M_\eq$. Let $\rho \geq 0$ be sufficiently large, so that the mapping $ q + {\rho}\|\cdot\|^2$ is convex. (Existence of such $\rho$ follows from \eqref{eq:q-weak-convex}.) Then, the function $f(y) := \frac12 q(Ty) + \frac{\rho}{2}\|Ty\|^2$ is also convex and using $\dom(q) \subset \eq$, $\mathrm{range}(P) = \eq$, the orthogonality of $P$ and $Q$ and $y = [ y_\eq^\top, \;\, y_\perp^\top ]^\top$, it follows
\[ f(y) = f(y_\eq,y_\perp) =  \frac12 q(Py_\eq + Qy_\perp) + \frac{\rho}{2}[\|y_\eq\|^2 + \|y_\perp\|^2] = \frac12 q(Py_\eq) + \iota_{\{0\}}(y_\perp) + \frac{\rho}{2}\|y_\eq\|^2. \]
Setting $f_\eq(y_\eq) := \frac12 q(Py_\eq) + \frac{\rho}{2}\|y_\eq\|^2 = \frac12 q_\eq(y_\eq) + \frac{\rho}{2}\|y_\eq\|^2$ and for $y_\eq \in \dom(q_\eq)$, this implies
\begin{align*} 
\begin{bmatrix} g_\eq \\ g_\perp \end{bmatrix} \in \partial f(y_\eq,0) \;\, & \iff \;\, f_\eq(y_\eq^\prime) + \iota_{\{0\}}(y_\perp^\prime) - f_\eq(y_\eq) \geq \iprod{g_\eq}{y_\eq^\prime - y_\eq} + \iprod{g_\perp}{y_\perp^\prime} \quad \forall~y_\eq^\prime \in \R^r, \; y_\perp^\prime \in \R^{n-r} \\ & \iff \;\, f_\eq(y_\eq^\prime) - f_\eq(y_\eq) \geq \iprod{g_\eq}{y_\eq^\prime - y_\eq} \quad \forall~y_\eq^\prime \in \R^r \;\, \iff \;\, g_\eq \in \partial f_\eq(y_\eq)
\end{align*}
and we have $\partial f(y_\eq,0) = \partial f_\eq(y_\eq) \times \R^{n-r}$. If $y_\eq \notin \dom(q_\eq)$, both subdifferentials are empty, so the identity remains valid. Hence, applying the change-of-coordinates and perturbation rules \cite[Exercises~10.7 and 10.10]{RockafellarWets}, we can infer 
\begin{equation} \label{eq:ortho-m} T^\top M(Py_\eq) = T^\top \partial (\tfrac12 q)(P y_\eq) = \partial(\tfrac12 q \circ T)(y_\eq,0) =  \partial(\tfrac12 q_\eq)(y_\eq) \times \R^{n-r} = M_\eq(y_\eq) \times \R^{n-r} \end{equation}
and $M^{-1}(u) = PM_\eq^{-1}(P^\top u)$, $u \in \R^n$. Here, notice that $M(z) = \varnothing$ if $Q^\top z = z_\perp \neq 0$, i.e., if $z \notin \eq$, and both sides of \eqref{eq:ortho-m} are empty when $y_\eq \notin \dom(q_\eq)$. As before, these identities imply that the metric and strong metric regularity of $M$ are equivalent to the corresponding properties of $M_\eq$. We record our observations in the following lemma. 

\begin{lemma}[Intrinsic reduction] \label[lemma]{lem:intrinsic-red}
    Let $q : \R^n \to \Rex$ be proper, lsc, and suppose that $q$ is $2\rho$-weakly convex for some $\rho \geq 0$. Set $\eq := \Span(\dom(q))$, $r := \dim(\eq)$, and let the columns of $P \in \R^{n \times r}$ form an orthonormal basis of $\eq$. Fix $\lambda > 0$ with $\lambda\rho < 1$ and define
\[ H := \nabla \env{q/2}{\lambda}, \quad M := \partial(\tfrac12 q), \quad q_\eq := q \circ P, \quad H_\eq := \nabla \env{q_\eq/2}{\lambda}, \quad M_\eq := \partial (\tfrac12 q_\eq). \]
    %
\begin{enumerate}[label=\textup{(\roman*)}]
\item $H$ is (strongly) metrically regular at $(\bar y,\bar u) \in \gph(H)$ iff $H_\eq$ is (strongly) metrically regular at $(P^\top\bar y,P^\top \bar u)$.
\item For $(\bar y,\bar u) \in \gph(M)$, we necessarily have $\bar y \in \eq$ and $M$ is (strongly) metrically regular at $(\bar y,\bar u)$ if and only if $M_\eq$ is (strongly) metrically regular at $(P^\top\bar y,P^\top \bar u)$.
\end{enumerate}
\end{lemma}

\Cref{lem:intrinsic-red} justifies replacing the models and functions $(q,H,M)$ by the intrinsic objects $(q_\eq,H_\eq,M_\eq)$. In particular, without loss of generality, we may fully work on $\eq$ and assume (following the orthogonal identifications in \eqref{eq:ortho-1}-\eqref{eq:ortho-2} and after relabeling of the dimension) that $\Span(\dom(q)) = \R^n$. We will leverage this reasoning more explicitly in \Cref{sec:euler,sec:lineality}. 

\section{An Euler characteristic-based criterion for strong metric regularity}
\label{sec:euler}

In this section, we formally introduce and prove the topological curvature criterion for strong metric regularity as stated in \Cref{thm:intro-euler}. Our discussion builds on the general setting and notation of the previous section, which we repeat and summarize here for ease of exposition 
%
\begin{subequations}\label{eq:euler-section-notation}
\begin{align}
q&:=\rd^2h(\bar x|\bar v),
\qquad H:= \nabla\env{q/2}{\lambda}, \qquad M := \partial(\tfrac12 q),
\label{eq:euler-section-models}\\
\eq &:=\Span(\dom(q)),
\qquad \Sph_\eq:= \mathbb S_{n-1} \cap \eq,
\label{eq:euler-section-space} \qquad \negq :=\{w\in\dom(q)\cap\Sph_\eq:q(w)<0\}.
\end{align}
\end{subequations}

If $\eq=\{0\}$, then $q=\iota_{\{0\}}$, $\env{q/2}{\lambda} = \norm{z}^2/(2\lambda)$, and $H(z)=\lambda^{-1}z$. In this case, all assertions below will hold with degree one and the sphere $\Sph_\eq$ and the negative link $\negq$ in \eqref{eq:euler-section-space} are empty. Thus, without loss of generality, we may assume $\dim(\eq)\ge1$.


\paragraph{Roadmap.} The derivation of the Euler characteristic-based criterion relies on three key arguments. 
First, we show that the set $\negq$ of negative curvature directions of $q$ is homotopy equivalent to the set $\mathscr M_q^\lambda := \{z \in \Sph_\eq: \env{q/2}{\lambda}(z) < 0\}$. Second, we use a homogeneous gradient-index formula to link the degree of $\nabla \env{q/2}{\lambda}$ to the Euler characteristic of $\mathscr M_q^\lambda$. Finally, the degree formula is combined with definable openness and discreteness: having the degree $\pm1$ is equivalent to a single local inverse branch, and metric regularity makes this branch Lipschitz continuous.


\subsection{Connecting ``negative'' directions}

As a first step, we establish that the sets of ``negative'' directions $\mathscr M_q^\lambda$ and $\negq$ have the same Euler characteristic.

\begin{lemma}[A proximal homotopy between $\mathscr M_q^\lambda$ and $\negq$]
\label[lemma]{lem:negative-link-homotopy}
%
Let $\sL$ be a finite-dimensional Euclidean space and let $q : \sL \to \Rex$ be a lsc, proper, definable, $2\rho$-weakly convex, $\rho \geq 0$, and positively $2$-homogeneous function. Set $\eq :=\Span(\dom(q))$ and $\Sph_\eq:=\{ x\in \eq:\norm{x}=1\}$. For fixed $\lambda > 0$ with $\lambda\rho < 1$, the sets
\[
\mathscr M_q^\lambda := \{z\in\Sph_\eq: \env{q/2}{\lambda}(z)<0\}
\quad\text{and}\quad
\negq := \{w\in\dom(q)\cap\Sph_\eq:q(w)<0\}
\]
are definable and have the same Euler characteristic.
\end{lemma}
\begin{proof} As a $2\rho$-weakly convex function, $q$ is prox-bounded and hence, definability of $\env{q/2}{\lambda}$, $\lambda\rho < 1$, follows from \Cref{lem:moreau-transform}~(iii). This immediately establishes definability of the sets $\mathscr M_q^\lambda$ and $\negq$. 
Applying \eqref{eq:cone-link-euler}, we further have $\chi(\mathscr M_q^\lambda) = \chi(\mathscr M^-)$ and $\chi(\negq) = \chi(\mathscr N^-)$, where
\[ \mathscr M^- := \{z \in \eq: \env{q/2}{\lambda}(z)<0\} \quad \text{and} \quad \mathscr N^- := \{w \in \eq: q(w) < 0\}. \]
Thus, it is sufficient to compare the cones $\mathscr M^-$ and $\mathscr N^-$.
%
Due to $\env{q/2}{\lambda}(z) \le q(z)/2$, we have $\mathscr N^-\subseteq \mathscr M^-$. Set $p:=\prox_{\lambda q/2}|_\eq: \eq \to \eq$. The weak convexity of $q$ ensures that the properties in \Cref{prop:prox-env}~(i) hold globally. In particular, $p$ is Lipschitz continuous and definable. 
%
For $z\in \mathscr M^-$, we obtain $\env{q/2}{\lambda}(z) = \half q(p(z))+\frac1{2\lambda}\norm{p(z)-z}^2 < 0$ and thus, $p(z)\in \mathscr N^-$.
%
We now define the mappings 
\[
f:\mathscr M^-\to \mathscr N^-, \quad f(z) := p|_{\mathscr M^-}(z),
\qquad
g:\mathscr N^-\to\mathscr M^-, \quad g(w) := w,
\]
and $\mathscr H:\mathscr M^-\times[0,1]\to \mathscr M^-$, $\mathscr H(z,t):=(1-t)z+tp(z)$. The function $\mathscr H$ is well-defined because, for $z\in \mathscr M^-$ and $t\in[0,1]$, we have
\begin{align*}
\env{q/2}{\lambda}(\mathscr H(z,t))
&\le \frac12 q(p(z))
  +\frac1{2\lambda}\norm{p(z)-\mathscr H(z,t)}^2 
=\frac12 q(p(z))
  +\frac{(1-t)^2}{2\lambda}\norm{p(z)-z}^2 \le \env{q/2}{\lambda}(z)<0.
\end{align*}
Moreover, it holds that $\mathscr H(\cdot,0)=\operatorname{id}_{\mathscr M^-}$ and $\mathscr H(\mathscr M^-,1)\subseteq \mathscr N^-$. Next, we verify that $\mathscr H(\mathscr N^-,t)\subseteq \mathscr N^-$ for all $t\in[0,1]$. For $z\in \mathscr N^-$, set $e_z(w) := \frac12 q(w) + \frac{1}{2\lambda} \|w-z\|^2$.
The function $e_z$ is (strongly) convex and has the unique minimizer $p(z)$. This yields 
\begin{align*}
e_z(\mathscr H(z,t))
&\le (1-t) e_z(z)+t e_z(p(z)) 
\le e_z(z)=\frac12 q(z)<0
\end{align*}
and $q(\mathscr H(z,t))<0$. Thus, it follows $\mathscr H(\mathscr N^-,t)\subseteq \mathscr N^-$ for every
$t\in[0,1]$.  We have therefore verified the following homotopy properties
\[
\mathscr H(\cdot,0)=\operatorname{id}_{\mathscr M^-},
\qquad
\mathscr H(\cdot,1)=g\circ f \quad \implies \quad g \circ f \simeq \operatorname{id}_{\mathscr M^-}
\]
and, for
$\mathscr H|_{\mathscr N^-\times[0,1]}:\mathscr N^-\times[0,1]\to \mathscr N^-$
\[
\mathscr H(\cdot,0)|_{\mathscr N^-}=\operatorname{id}_{\mathscr N^-},
\qquad
\mathscr H(\cdot,1)|_{\mathscr N^-}=f\circ g \quad \implies \quad f \circ g \simeq \operatorname{id}_{\mathscr N^-}.
\]
Hence, the sets $\mathscr M^-$ and $\mathscr N^-$ are homotopy equivalent and we have $\chi(\mathscr N^-) =\chi(\mathscr M^-)$. This concludes the proof.
\end{proof}

\subsection{Gradient degree formula}

The next lemma provides a variant of the Khimshiashvili formula recalled in
\cite[Theorem~2.3 and Corollary~2.4]{Dutertre2020} for positively two-homogeneous, definable $C^1$-functions. The result relates the local Brouwer degree of a gradient mapping to the Euler characteristic of its negative directions. 

\begin{lemma}[Gradient-degree formula for positively $2$-homogeneous functions]
\label[lemma]{lem:homogeneous-index}
Let $\sL$ be a (finite-dimensional) Euclidean space and let $f: \sL\to\R$ be definable, $C^1$, and positively $2$-homogeneous.  If
$0$ is isolated in $(\nabla f)^{-1}(0)$, then
\begin{equation}
\deg_0(\nabla f)
=1-\chi(\{z\in\sphL: f(z)<0\}).
\label{eq:khimshiashvili-homogeneous}
\end{equation}
\end{lemma}

Dutertre states and proves the formula \eqref{eq:khimshiashvili-homogeneous} for definable $C^r$-functions with $r\ge2$ \cite{Dutertre2020}. Here, we extend this setting to positively two-homogeneous, definable $C^1$-functions. This generalization is necessary in our context, since we plan to apply \cref{lem:homogeneous-index} to the Moreau envelope $\env{q/2}{\lambda}$ which is only continuously differentiable. 

In the subsequent proof, we approximate the restriction of $f$ to the sphere by a definable $C^2$-function and verify that this perturbation preserves both the local degree and the Euler characteristic of the negative set. As the proof of \Cref{lem:homogeneous-index} is somewhat involved, we have decided to verify several more technical steps in a separate section in the appendix. Upon first reading of the manuscript, the reader may safely skip to \Cref{sec:main-results}.

\begin{proof}
In the following, for a function $h : \sphL \to \R$, we will use the abbreviation $\{h\,\circ\, 0\} := \{z \in \sphL : h(z)\,\circ\, 0\}$ where $\circ \in \{\leq,<,=\}$.
Positive $2$-homogeneity and differentiability at the origin give $f(0)=0$ and $\nabla f(0)=0$ by \cite[Section~1]{Haraux2005}. If $\sL=\{0\}$, the local degree is one by the zero-dimensional convention, while $\sphL=\varnothing$ and the right-hand side of \eqref{eq:khimshiashvili-homogeneous} is also one. Set $d:=\dim(\sL)\ge1$, so $\dim(\sphL)=d-1$, and let $\phi:=f|_{\sphL}$.  For $u\in\sphL$, the tangent space $\mathcal T_u\sphL$ and the orthogonal projection onto $\mathcal T_u\sphL$ are given by
\[ \mathcal T_u\sphL=\{z\in \sL:\iprod{z}{u}=0\}, \qquad \mathcal P_{\mathcal T_u\sphL}=I-uu^\top.
\]
Thus, the Riemannian gradient of $\phi$ on $\sphL$ is
\[ \nabla_{\sphL}\phi(u):=\mathcal P_{\mathcal T_u\sphL}\nabla f(u)=(I-uu^\top)\nabla f(u), \qquad u \in \sphL.
\]
The differential $\rd \phi(u) : \mathcal T_u\sphL \to \R$ is then characterized through the identity $\rd \phi(u)[\xi] = \iprod{\nabla_{\sphL}\phi(u)}{\xi} = \iprod{\nabla f(u)}{\xi}$ for all $\xi \in \mathcal T_u\sphL$, \cite{absil2008manifold}. Euler's identity for the positively two-homogeneous function $f$ gives
\begin{equation}
\iprod{\nabla f(u)}{u}=2f(u)=2\phi(u).
\label{eq:euler-identity-phi}
\end{equation}
Decomposing $\nabla f(u)$ into its radial and tangential parts, we then obtain
\begin{equation}
\begin{aligned}
\nabla f(u)
&=\iprod{\nabla f(u)}{u}u+\nabla_{\sphL}\phi(u) =2\phi(u)u+\nabla_{\sphL}\phi(u).
\end{aligned}
\label{eq:radial-tangential-gradient}
\end{equation}
Positive two-homogeneity also yields $\nabla f(tu)=t\nabla f(u)$ for $t>0$. Hence, the right-hand side of \eqref{eq:radial-tangential-gradient} is nonzero: Otherwise, the ray through $u$ would contain nonzero stationary points arbitrarily close to the isolated stationary point $0$.
Thus, $0$ is a regular value of $\phi$, i.e., every $u\in\sphL$ with $\phi(u)=0$ satisfies $\rd\phi(u)\ne0$ or, equivalently, $\nabla_{\sphL}\phi(u)\ne0$.

We first consider the case $d=1$. Let us write $\sphL=\{e,-e\}$ and set $a_+:=\phi(e)$ and $a_-:=\phi(-e)$. Since the tangential gradients vanish, \eqref{eq:radial-tangential-gradient} gives $\nabla f(e)=2a_+e$ and $\nabla f(-e)=-2a_-e$; in particular, $a_+a_-\ne0$. By the one-dimensional local degree formula \cite[Proposition~1.2.3]{DincaMawhin}, we then have
\[
\deg_0(\nabla f)
=\frac{\sgn(a_+) + \sgn(a_-)}{2}
=1-|\{\sigma\in\{+,-\}:a_\sigma<0\}|.
\]
The last cardinality is exactly the Euler characteristic of the set $\{z \in \{e,-e\}: f(z) < 0\}$. Hence, \eqref{eq:khimshiashvili-homogeneous} holds if $d=1$, and we may assume that $d\ge2$.

For $h\in C^1(\sphL)$, we introduce the norm $\norm{h}_{C^1(\sphL)} :=\norm{h}_\infty+\norm{\nabla_{\sphL}h}_\infty$ and define
\[
V_h(u):=2h(u)u+\nabla_{\sphL}h(u).
\]
By \eqref{eq:radial-tangential-gradient}, we have $V_\phi=\nabla f$ on $\sphL$, so $V_\phi$ is continuous and nowhere zero. Thus, since $\sphL$ is compact, it holds that
\[
\delta:={\min}_{u\in\sphL}~\norm{V_\phi(u)}>0.
\]
The definable approximation theorem \cite[Theorem~1.1]{Escribano2002} gives a definable $C^2$-function $\widetilde\phi:\sphL\to\R$ arbitrarily close to $\phi$ in the $\norm{\cdot}_{C^1(\sphL)}$-norm. Choose $\widetilde\phi$, so that
\begin{equation}
2\norm{\widetilde\phi-\phi}_\infty
+\norm{\nabla_{\sphL}(\widetilde\phi-\phi)}_\infty
<\frac{\delta}{2}.
\label{eq:c2-approximation-size}
\end{equation}
For $\phi_s:=(1-s)\phi+s\widetilde\phi$, $0\le s\le1$, and using the definition of $h\mapsto V_h$ and \eqref{eq:c2-approximation-size}, it follows
\begin{align*}
\norm{V_{\phi_s}(u)-V_\phi(u)}
&=s\norm{2(\widetilde\phi(u)-\phi(u))u
 +\nabla_{\sphL}(\widetilde\phi-\phi)(u)} \le 2\norm{\widetilde\phi-\phi}_\infty
 +\norm{\nabla_{\sphL}(\widetilde\phi-\phi)}_\infty
<\frac{\delta}{2}.
\end{align*}
This implies
\begin{equation}
\norm{V_{\phi_s}(u)}
\ge\norm{V_\phi(u)}-\norm{V_{\phi_s}(u)-V_\phi(u)}
>\frac{\delta}{2}
\label{eq:homotopy-field-lower-bound}
\end{equation}
for every $(u,s)\in\sphL\times[0,1]$.  In particular, we have
\begin{equation}
\norm{\nabla_{\sphL}\phi_s(u)}
=\norm{V_{\phi_s}(u)}
>\frac{\delta}{2}
\quad \forall~u\in\sphL,\ s\in[0,1],\ \quad \text{with} \quad \phi_s(u)=0.
\label{eq:zero-level-gradient-bound}
\end{equation}
Thus, $0$ is a regular value of $\phi_s$ for all $s\in[0,1]$. Set $B := \{z \in \sL: \|z\| < 1\}$. For $s\in[0,1]$, we define the conical extension $\mathcal V_s:\cl(B) \to \sL$ of $V_{\phi_s}$ by
\[
    \mathcal V_s(z):=
    \begin{cases}
        \|z\| \cdot V_{\phi_s}(z/\|z\|),&z\ne0,\\
        0,&z=0.
    \end{cases}
\]
The map $(z,s)\mapsto \mathcal V_s(z)$ is continuous because $V_{\phi_s}$ is uniformly bounded on $\sphL\times[0,1]$. By positive homogeneity, we have $\mathcal V_0(z) = \|z\| V_\phi(z/\|z\|) = \|z\|\nabla f(z/\|z\|) = \nabla f(z)$ for all $z$. Moreover, \eqref{eq:homotopy-field-lower-bound} implies $\mathcal V_s(u)\ne0$ for all $u\in \sphL = \partial B$ and $s\in[0,1]$. Hence, homotopy invariance of the Brouwer degree (cf.\ \ref{D3}) yields
\begin{equation}\label{eq:degree-approximation}
    \deg_0(\nabla f) = \deg(\mathcal V_0,B,0) = \deg(\mathcal V_1,B,0).
\end{equation}
Next, we consider the definable function $\Psi:\sL\to\R$ defined by 
\[
\Psi(z) :=
\begin{cases}
\norm{z}^4\cdot \widetilde\phi(z/\norm{z}),&z\ne0,\\
0,&z=0.
\end{cases}
\]
The function $\Psi$ is positively four-homogeneous and $C^2$ on $\sL\setminus\{0\}$. For $k=0,1,2$, differentiating this homogeneity identity shows that $D^k\Psi$ is positively $(4-k)$-homogeneous (on $\sL\setminus\{0\}$). Its restriction to the compact sphere $\sphL$ is bounded, so there exists $C_0>0$ such that
\[
\norm{D^k\Psi(z)}\le C_0 \norm{z}^{4-k}
\qquad \forall~z \ne0 \quad \text{and} \quad k=0,1,2.
\]
The estimates for $k=0$ and $k=1$ show successively that
$D\Psi(0)=0$ and $D^2\Psi(0)=0$; the estimate for $k=2$ makes the Hessian
continuous at $0$ and hence, it follows $\Psi\in C^2(\sL)$.  Moreover, for $r>0$ and $u\in\sphL$, we obtain
\begin{equation}
\nabla\Psi(ru)
=r^3 (4\widetilde\phi(u)u
          +\nabla_{\sphL}\widetilde\phi(u)).
\label{eq:quartic-gradient}
\end{equation}
For $u\in\sphL$, the vectors $u$ and $\nabla_{\sphL}\widetilde\phi(u)$ are orthogonal. Consequently, for all $(u,t)\in\sphL\times[2,4]$, it holds that
\begin{equation} \label{eq:another-homotopy}
\norm{t\widetilde\phi(u)u
      +\nabla_{\sphL}\widetilde\phi(u)}^2
=t^2\widetilde\phi(u)^2
 +\norm{\nabla_{\sphL}\widetilde\phi(u)}^2>0.
\end{equation}
Indeed, equality would imply both $\widetilde\phi(u)=0$ and $\nabla_{\sphL}\widetilde\phi(u)=0$, which implies $V_{\widetilde \phi}(u) = V_{\phi_1}(u) = 0$ and contradicts \eqref{eq:homotopy-field-lower-bound}.  Taking $t=4$ in \eqref{eq:another-homotopy} and using \eqref{eq:quartic-gradient} shows that $0$ is an isolated zero of $\nabla\Psi$. 
Consider the continuous homotopy
\[ \mathscr H(z,\tau) := (1-\tau)\mathcal V_1(z)+\tau\nabla\Psi(z),
    \qquad (z,\tau)\in \cl(B)\times[0,1].
\]
For $u\in \sphL = \partial B$, we have $\mathscr H(u,\tau) =(2+2\tau)\widetilde\phi(u)u +\nabla_{\sphL}\widetilde\phi(u)$.
Thus, invoking \eqref{eq:another-homotopy}, we can infer $\mathscr H(u,\tau) \neq 0$ for all $(u,\tau) \in \partial B \times [0,1]$. Homotopy invariance therefore gives
\[
    \deg(\mathcal V_1,B,0)
    =
    \deg(\nabla\Psi,B,0)
    =
    \deg_0(\nabla\Psi).
\]
The spherical-sublevel formula \cite[Corollary~2.4]{Dutertre2020}, applied to the $C^2$-function $\Psi$ (and using the positive homogeneity of $\Psi$), yields
\[
    \deg_0(\nabla\Psi) = 1 - \chi(\{z \in \sphL: \Psi(z) \leq 0\})
    = 1-\chi(\{\widetilde\phi\le0\}).
\]
Since $0$ is a regular value of $\widetilde\phi$, \cite[Proposition~2.1]{Jubin2023} shows that $\mathcal C:=\{\widetilde\phi\le0\}$ is a full-dimensional $C^2$-submanifold with boundary of $\sphL$, and that its manifold boundary, which agrees with its topological boundary in $\sphL$, is $\partial \mathcal C=\{\widetilde\phi=0\}$. Consequently, the interior of $\mathcal C$ is $\mathcal C\setminus\partial \mathcal C=\{\widetilde\phi<0\}$.  A topological manifold with boundary is homotopy equivalent to its interior \cite[Proposition~3.11]{GalashinKarpLam2022}.  Hence the homotopy invariance of the Euler characteristic gives
\[
\chi(\{\widetilde\phi\le0\})
=\chi(\{\widetilde\phi<0\}).
\]
So far, we have shown that $\deg_0(\nabla f) = \deg_0(\nabla\Psi) = 1-\chi(\{\widetilde\phi<0\})$. Thus, the proof is complete if we can establish $\chi(\{\widetilde\phi<0\}) = \chi(\{\phi<0\})$. For ease of exposition, as the verification of this identity is rather technical and somewhat independent of the previous steps, we defer this remaining part of the proof to \Cref{sec:app:remainder-of-lemma}.
\end{proof}

\subsection{Main result} \label{sec:main-results}

We now present and prove the main result of this section. 

\begin{theorem}[An Euler-link curvature criterion for strong metric regularity]
\label{thm:abstract-euler}
Let $h:\R^n\to\Rex$ be proper, lsc, and definable. Let $(\bar x,\bar v)\in\gph(\partial h)$ be given and assume that $h$ is prox-regular and subdifferentially continuous at $\bar x$ for $\bar v$.  If $\partial h$
is metrically regular at $(\bar x,\bar v)$, then
\begin{equation}
\partial h\text{ is strongly metrically regular at }(\bar x,\bar v)
\quad\iff\quad
\chi(\negq)\in\{0,2\}.
\label{eq:main-euler-box}
\end{equation}
\end{theorem}

Before verifying \Cref{thm:abstract-euler}, we provide an additional computational result that links the local degree of $H$ and the Euler characteristic of the set $\negq$. This connection is a consequence of the gradient index formula.


\begin{proposition}
\label[proposition]{prop:model-degree}
Let $h$ satisfy the conditions stated in \Cref{thm:abstract-euler}. For the model $H$ in
\eqref{eq:euler-section-models}, we then have
\begin{equation}
\deg_0(H) = \deg_0(\nabla \env{q/2}{\lambda}) =1-\chi(\negq).
\label{eq:model-degree}
\end{equation}
\end{proposition}

\begin{proof}
By the metric regularity of $H$ at $(0,0)$ in \Cref{thm:exact-reduction} and the equivalence between metric regularity and linear openness (see \Cref{sec:intro-metric}), after shrinking to a connected definable open neighborhood $V$ of the origin, the restriction $H|_V$ is open.  Since $H|_V$ is continuous and definable, \Cref{thm:definable-openness} shows that it has finite fibers. Since $H^{-1}(0)$ is a cone, this finiteness forces
$H^{-1}(0)=\{0\}$ near the origin. Using the orthogonal splitting
\eqref{eq:ortho-1}--\eqref{eq:ortho-2} and \eqref{eq:deg-ortho}, we have $\deg_0(H)=\deg_0(H_\eq)$. Applying \Cref{lem:homogeneous-index} to the function $e_\eq = \env{q_\eq/2}{\lambda}$ on $\sL = \R^r$, $r = \dim(\eq)$, noticing $\Span(\dom(q_\eq)) = \R^r$, and using \Cref{lem:negative-link-homotopy}, we can infer 
\[ \deg_0(H_\eq) = \deg_0(\nabla e_\eq) = 1-\chi(\{y \in \Sph_{r-1}: e_\eq(y) <0\}) = 1- \chi(\mathscr N_{q_\eq}),
\] 
where $\mathscr N_{q_\eq} = \{y \in \dom(q_\eq) \cap \Sph_{r-1}: q_\eq(y) < 0 \}$. Finally, let $P \in \R^{n \times r}$ be the orthonormal basis of $\eq$---as introduced in \Cref{sec:regularity-moreau} and define $f : \mathscr N_{q_\eq} \to \negq$ by $f(y) = Py$ and $g : \negq \to \mathscr N_{q_\eq}$ by $g(z) = P^\top z$. These maps are well-defined and continuous: $P$ is a linear isometry from $\R^r$ to $\eq$ and we have $q_\eq(y) = q(Py)$. Moreover, it holds that $(g\circ f)(y) = P^\top P y = y$ and $(f\circ g)(z) = PP^\top z$ for $y \in \R^r$ and $z \in \eq$. This shows $g \circ f = \mathrm{id}_{\mathscr N_{q_\eq}}$ and $f \circ g = \mathrm{id}_{\negq}$ and consequently, the sets $\mathscr N_{q_\eq}$ and $\negq$ are homotopy equivalent. Combining these steps, we therefore obtain $\deg_0(H) = 1 - \chi(\negq)$---as desired.
\end{proof}

We can now turn to the proof of \Cref{thm:abstract-euler}.

\begin{proof}[Proof of \Cref{thm:abstract-euler}]
By \Cref{thm:exact-reduction}, strong metric regularity of $\partial h$ is equivalent to strong metric regularity of the map $H = \nabla \env{q/2}{\lambda}$. As shown in the proof of \Cref{prop:model-degree}, $H$ is open and discrete on a neighborhood of the origin. Hence, the homological inverse-mapping criterion, cf.\ \Cref{lem:open-degree-one}, together with the fact that a local homeomorphism has local degree $\pm1$, see \cite[Section~2.2, p.~136]{Hatcher}, gives
\[
H\text{ is strongly metrically regular at }(0,0)
\quad\Longleftrightarrow\quad
|\deg_0(H)|=1.
\]
Indeed, if $|\deg_0(H)|=1$, then \Cref{lem:open-degree-one} implies that the map $H$ is a local homeomorphism at the origin. Since $H$ is metrically regular at $(0,0)$, its local inverse is Lipschitz. Finally, by \eqref{eq:model-degree}, we have $|\deg_0(H)|=1$ if and only if $|1-\chi(\negq)|=1$, i.e., $\chi(\negq)\in\{0,2\}$. This completes the proof.
\end{proof}

Next, we record a formal result for the observations made in \eqref{eq:intro-second-order-stationarity}--\eqref{eq:intro-second-order-consequence}. If $\bar x$ is a second-order stationary point of the tilted function
$h_{\bar v} = h-\langle\bar v,\cdot\rangle$, then $\negq=\varnothing$ and $\partial h$ is strongly metrically regular at $(\bar x, \bar v)$. 

\begin{corollary}[Second-order stationarity and strong metric regularity]
\label[corollary]{cor:second-order-stationarity}
Let the assumptions stated in \Cref{thm:abstract-euler} hold and suppose that  $q(w) = \mathrm{d}^2h(\bar x|\bar v)(w)\ge0$ for all $w\in\dom(q)$, i.e., $\bar x$ is a second-order stationary point of $h_{\bar v}$. Then, $\partial h$ is strongly metrically regular at $(\bar x,\bar v)$.
\end{corollary}

We conclude this section and present two examples to illustrate how the curvature criterion can be applied in practice.

\begin{example}[Generalized quadratic second subderivatives]
We first revisit \Cref{ex:structured-models} and assume that the second subderivative $q$ has the generalized quadratic form $q(w)=\langle w,Qw\rangle+\iota_K(w)$ as in \eqref{eq:example-quadratic}. Here, $Q\in\mathbb S^n$ and $K$ is a linear subspace. Let the columns of $P\in\mathbb R^{n\times r}$ form an orthonormal basis of $K$, and diagonalize $A:=P^\top QP=U\Gamma U^\top$, $\Gamma=\operatorname{diag} (\gamma_1,\ldots,\gamma_r)$, with $U \in \R^{r \times r}$ orthogonal. Suppose that the eigenvalues of $A$ are ordered so that $\gamma_1,\ldots,\gamma_m<0$ and $\gamma_i \geq 0$, $i = m+1,\dots,r$, where $m:=n_-(A)$. Then, the sets $\negq$ and
\[
    \mathscr N_\Gamma
    :=\{y\in\mathbb S_{r-1}: y^\top\Gamma y<0\}.
\]
are homotopy equivalent. (This can be shown following the proof of \Cref{prop:model-degree} and using the isometry $y\mapsto PUy$, with inverse $w\mapsto U^\top P^\top w$.) If $m=0$, then $\mathscr N_\Gamma=\varnothing$.
For $m\geq1$, set $\mathscr T:= \mathbb S_{m-1}\times\{0\} \subseteq\mathbb R^m\times\mathbb R^{r-m}$.
Since the first $m$ eigenvalues of $A$ are negative, we have $\mathscr T\subseteq\mathscr N_\Gamma$. Moreover, every $y=(y_-,y_{\geq0})\in\mathscr N_\Gamma$ satisfies $y_-\neq0$. Thus, the maps $f:\mathscr N_\Gamma\to \mathscr T$, $f(y):=({y_-}/{\|y_-\|},0)$ and $g: \mathscr T\to\mathscr N_\Gamma$, $g(u):=u$ are well-defined and continuous, with
$f\circ g=\mathrm{id}_{\mathscr T}$. Define
\[
    \mathscr H:\mathscr N_\Gamma\times[0,1]\to\mathscr N_\Gamma,
    \qquad \mathscr H(y,t):= \frac{(y_-,(1-t)y_{\geq0})}{\|(y_-,(1-t)y_{\geq0})\|}.
\]
The denominator is nonzero because $y_-\neq0$ and the normalization preserves strict negativity. Hence, $\mathscr H$ is well-defined and continuous and satisfies $\mathscr H(y,0)=y$, $\mathscr H(y,1)=(g\circ f)(y)$, and $\mathscr H(u,t)=u$ for every $u\in \mathscr T$.
Consequently, $g\circ f\simeq\mathrm{id}_{\mathscr N_\Gamma}$ and the sets $\negq$, $\mathscr N_\Gamma$, and $\mathscr T$ are all homotopy equivalent. The standard Euler characteristic of a sphere and the product formula \cite[Section~2.2]{Hatcher} then give, (including $m=0$),
\[
    \chi(\negq)=1-(-1)^m\in\{0,2\}.
\]
Thus, under the hypotheses of Theorem~\ref{thm:abstract-euler}, metric regularity of $\partial h$ always implies strong metric regularity in this case. 
\end{example}

\begin{example}[Failure of strong metric regularity]
\label{rem:ma-sun-euler-check}
We now discuss a counterexample by Ma and Sun presented in \cite[Section~5]{MaSun}. Define
\[ f: \R^2 \to \R, \quad f(x) = f(r(\cos(\theta),\sin(\theta))):=r^2\cos(3\theta),
\qquad
g : \R^2 \to \R, \quad g(x):=4\norm{x}^2+f(x),
\]
set $S:=\epi(g)$, and consider
\[
h(x,t):=-4\norm{x}^2+\iota_S(x,t),
\qquad
\bar v:=(0,0,-1).
\]
Ma and Sun show that the subdifferential $\partial h$ is metrically regular at $(0,\bar v)$ but has two local solution branches, and hence it is not strongly metrically regular. Here, we will find that $\negq$ has three connected components, so $\chi(\negq)=3$ and by \eqref{eq:model-degree}, it follows $\deg_0(H)=-2$. Thus, the topology of the second-order model records the failure of a unique inverse branch. We first compute $q = \mathrm d^2 h(0|\bar v)$. Since $g$ is continuous, $S$ is closed. The indicator-function formula \cite[equation~(3.1)]{BenkoMehlitz2023} then gives
\begin{align}
{\mathrm d}^2\iota_S(0|\bar v)(w,\omega)
&=
\liminf_{\substack{\tau\downarrow0,\ (w_\tau,\omega_\tau)\to (w,\omega)\\
                    (\tau w_\tau,\tau \omega_\tau)\in S}}
-\frac{2\iprod{\bar v}{(w_\tau,\omega_\tau)}}{\tau}=
\liminf_{\substack{\tau\downarrow0,\ (w_\tau,\omega_\tau)\to(w,\omega)\\
                    \omega_\tau\geq\tau g(w_\tau)}}
\frac{2\omega_\tau}{\tau},
\label{eq:ma-sun-indicator-second-subderivative}
\end{align}
where the second equality uses $\bar v=(0,0,-1)$ and $g(\tau w_\tau)=\tau^2g(w_\tau)$. Since $g(x)\geq3\norm{x}^2$, there is no feasible sequence when $\omega<0$, while for $\omega>0$ and every feasible sequence, we obtain $2\omega_\tau/\tau \to +\infty$. When $\omega=0$, feasibility gives $2\omega_\tau/\tau\geq2g(w_\tau)\to2g(w)$, and equality is achieved by taking $w_\tau=w$ and $\omega_\tau=\tau g(w)$. Therefore, we have
\[
{\mathrm d}^2\iota_S(0|\bar v)(w,\omega) = 2g(w)+\iota_{\{0\}}(\omega).
\]
The smooth-addition formula \cite[Exercise~13.18]{RockafellarWets}, applied to $(x,t)\mapsto-4\norm{x}^2$, now yields
\[ q(w,\omega) = {\mathrm d}^2h(0|\bar v)(w,\omega) =-8\norm{w}^2+2g(w)+\iota_{\{0\}}(\omega) = 2f(w)+\iota_{\{0\}}(\omega).
\]
On the set $\dom(q) \cap \Sph_\eq = \Sph_1 \times \{0\}$, negativity of $q$ is equivalent to the condition $\cos(3\theta)<0$. Therefore, the negative link $\negq$ can be represented as the following disjoint union 
\[ \negq = \left\{ x = (\cos(\theta), \sin(\theta))^\top: \theta \in \Big(\frac{\pi}{6},\frac{\pi}{2}\Big)
\mathbin{\dot\cup}
\Big(\frac{5\pi}{6},\frac{7\pi}{6}\Big)
\mathbin{\dot\cup}
\Big(\frac{3\pi}{2},\frac{11\pi}{6}\Big) \right\} \times \{0\},
\]
where each interval parameterizes an open arc. Each arc is contractible, so $H_0(\negq;\mathbb Q)\cong \mathbb Q^3$ and $H_j(\negq;\mathbb Q)=0$ for $j\geq1$. Hence, we obtain $\chi(\negq)=1+1+1=3$ and $\deg_0(H)=-2$ (by \eqref{eq:model-degree}). Consequently, \Cref{thm:abstract-euler} precisely recovers (and characterizes) the failure of strong metric regularity. More generally, Ma and Sun replace $f$ by $f_n(r(\cos(\theta),\sin(\theta)))=r^2\cos((n+1)\theta)$ and the coefficient $4$ in $g$ and $h$ by a sufficiently large constant $c_n$ \cite[Remark~5.4]{MaSun}. The resulting negative link has $n+1$ contractible components, and \eqref{eq:model-degree} gives degree
$-n$, which is in accordance with the observations reported in \cite{MaSun}.  


\end{example}

\section{Weak second-order stationarity and strong metric regularity}
\label{sec:lineality}

In this section, we show that for a general weakly convex and positively two-homogeneous model $q$, nonnegativity on the lineality space $\dom(q)\cap(-\dom(q))$ excludes branching under metric regularity and implies quadratic growth on the whole effective domain $\dom(q)$. This formalizes the results sketched in \Cref{thm:intro-lineality}.

\subsection{Uniform propagation of negative inertia}

The following lemma shows that the negative inertia of the Hessian propagates under metric regularity.

\begin{lemma} \label[lemma]{lem:uniform-mr-gradient-inertia}
Let $m\geq2$ and let $F=\nabla f:\R^m\to\R^m$ be Lipschitz, where $f : \R^m \to \R$ is definable and $C^1$. Suppose that $F$ is metrically regular at
$(x,F(x))$.  Then there exist $\epsilon > 0$ and $\kappa>0$ such that, upon
setting
\[
 \mathcal B_{F,\epsilon}(x) := \{z\in\bB_\epsilon(x):
    F\text{ is }C^1\text{ on a neighborhood of $z$} \},
\]
the negative inertia $z \mapsto \nminus(DF(z))$ is constant on $\mathcal B_{F,\epsilon}(x)$.  In addition, $\mathcal B_{F,\epsilon}(x)$ is definable and open and dense in $\bB_\epsilon(x)$ and we have
\begin{equation}
 \sigma_{\min}(DF(z))\geq\kappa^{-1}
 \qquad \forall~z\in\mathcal B_{F,\epsilon}(x).
 \label{eq:uniform-conorm}
\end{equation}
\end{lemma}

\begin{proof}
By the definition of metric regularity (cf.\ \eqref{eq:def-mr}) and the continuity of $F$, there exist $\epsilon>0$ and $\kappa>0$ such that
\begin{equation} \label{eq:step-mr-01}
 \dist(z,F^{-1}(y))
 \leq\kappa\norm{F(z)-y}
 \qquad \forall~z\in\bB_\epsilon(x), \quad y\in F(\bB_\epsilon(x))+\bB_\epsilon.
\end{equation}
For every $z\in\bB_\epsilon(x)$, the estimate \eqref{eq:step-mr-01} makes $F$ metrically regular at $(z,F(z))$ with modulus $\kappa$. Since $F$ is $C^1$ around each $z\in\mathcal B_{F,\epsilon}(x)$, we can infer that $DF(z)$ is nonsingular, and using \cite[Theorem~5D.5 and Corollary~5A.2]{dontchev2009implicit}, it follows
\[ \norm{DF(z)^{-1}}\leq\kappa. \]
Equivalently, we have $\sigma_{\min}(DF(z))\geq\kappa^{-1}$ for all $z\in\mathcal B_{F,\epsilon}(x)$, which proves \eqref{eq:uniform-conorm}.
Since $f$ is definable and $C^1$, $F = \nabla f$ is definable. By the equivalence between metric regularity and linear openness recalled in \Cref{sec:intro-metric}, \eqref{eq:step-mr-01} implies that $F|_{\bB_\epsilon(x)}$ is an open (and definable) map. Since the set $\bB_\epsilon(x)$ is connected, \cref{thm:definable-openness} is applicable and we can infer that the fibers of $F|_{\bB_\epsilon(x)}$ are finite and $z \mapsto \sgn(\det(DF(z)))$ is constant on $\mathcal B_{F,\epsilon}(x)$. By \cite[Lemma~3.3~(ii)]{DinhPham}, the set $\mathcal B_{F,\epsilon}(x)$ is open and dense in $\bB_\epsilon(x)$. Furthermore, combining \cite[Lemmas~2.7 and 2.5~(i)]{DinhPham}, $\mathcal B_{F,\epsilon}(x)$ is definable and has finitely many connected components.
Since $F=\nabla f$, the matrices $DF(z)$ are symmetric. Hence, since $DF$ is continuous on $\mathcal B_{F,\epsilon}(x)$, its negative inertia is constant on each connected component of $\mathcal B_{F,\epsilon}(x)$. 

We now show that these constants agree. To do so, we leverage a technique used by Dinh and Pham (in \cite{DinhPham}) to prove the definable open-mapping theorem. Let $\Gamma$ be the graph whose vertices are the connected components of $\mathcal B_{F,\epsilon}(x)$, and join two distinct components $\mathcal R_i$ and $\mathcal R_j$ by an edge if and only if
\begin{equation} \label{eq:rule-adjacency}
 \dim(\bB_\epsilon(x)\cap\cl({\mathcal R_i})\cap\cl({\mathcal R_j}))\geq m-1,
\end{equation}
where $\dim$ denotes the o-minimal dimension.\footnote{For a nonempty definable set $X$, the dimension $\dim(X)$ is the largest dimension of a manifold in a finite definable $C^1$-manifold decomposition of $X$. This value is independent of the chosen decomposition. By convention, $\dim(\varnothing) =-\infty$; see \cite[Section~2.3]{DinhPham}.} 
Then $\Gamma$ is connected by \cite[Fact~2.2]{PeterzilStarchenko2007}; see also \cite[Lemma~3.1]{DinhPham}.  
Let $\mathcal R_+$ and $\mathcal R_-$ be two adjacent components, and set
\[
 S:=\bB_\epsilon(x)\cap\cl({\mathcal R_+})\cap\cl({\mathcal R_-}).
\]
Since $\mathcal R_+$ and $\mathcal R_-$ are disjoint open connected sets, $S$ is contained in $\partial \mathcal R_+\cap\partial \mathcal R_-$. Moreover, since $\mathcal R_+$ is open, we have $\dim(\mathcal R_+) =m$, and applying \cite[Lemma~2.6~(iii)]{DinhPham}, it follows $\dim(S)<m$. The rule of adjacency \eqref{eq:rule-adjacency} gives $\dim(S)\geq m-1$, which implies
$\dim(S)=m-1$. Hence, a definable $C^1$-manifold decomposition of $S$ contains an $(m-1)$-dimensional embedded manifold $\mathcal M$ of $\R^m$.
In the following, we ``flatten'' the set $\mathcal M$ to facilitate the application of the tangential derivative-limit lemma provided in \cite{PeterzilStarchenko2007}. Choose $z_0\in\mathcal M$. By the implicit function theorem, there exist an open neighborhood $U \subset \bB_\epsilon(x)$ of $z_0$, an orthogonal matrix $Q$, a nonempty definable open set $W\subset\mathbb{R}^{m-1}$, and a definable $C^1$-function $\varphi:W\to\mathbb{R}$ such that
\[
\mathcal M \cap U
=
\{z_0+Q(y,\varphi(y)):y\in W\}.
\]
Replacing $\mathcal M$ by $\mathcal M\cap U$, we may therefore assume that $\mathcal M$ has such representation. Choosing $W$ and $\delta>0$ sufficiently small, the map
\[
\Phi:W\times(-\delta,\delta)\to N\subset \bB_\epsilon(x),
\qquad \Phi(y,t):=z_0+Q(y,\varphi(y)+t),
\]
is a definable $C^1$-diffeomorphism onto an open set $N$, with $\Phi(W\times\{0\})=\mathcal M$. A similar construction is also used in \cite[Theorem~3.2]{PeterzilStarchenko2007} and \cite[Theorem~3.1, (i)$\implies$(ii)]{DinhPham}. Next, let us define
\[
\mathcal U_\pm:=\Phi^{-1}(\mathcal R_\pm\cap N), \quad \mathcal U = \mathcal U_+\cup \mathcal U_- , \quad \widehat F:=F\circ\Phi,  \quad g(y):=F(\Phi(y,0)).
\]
Then $\mathcal U_\pm$ and $\mathcal U$ are definable open sets, $W\times\{0\}\subset\partial \mathcal U_+ \cap \partial\mathcal U_-$, and each $\widehat F_i$ is continuous on $W\times(-\delta,\delta)$ (hence on $\mathcal U \cup (W \times \{0\})$) and is $C^1$ on $\mathcal U$. The mapping $g:W\to\R^m$ is definable and without loss of generality, we may assume that $g$ is $C^1$ on $W$.\footnote{Applying \cite[Lemma~2.7]{DinhPham} with $p=1$, the set $D_g:= \{ a\in W: g \text{ is not } C^1 \text{ on any open neighborhood of } a \text{ contained in } W \}$ is definable with $\dim(D_g) < m-1$. Thus, $W\setminus D_g = \bigcup \{O \subset W: O \text{ is open and } g|_O \text{ is } C^1\}$ is nonempty and open and we can choose a nonempty definable open ball $W_0\subset W\setminus D_g$. Then $g|_{W_0}$ is $C^1$. Hence, restricting $\Phi$ to $W_0\times(-\delta,\delta)$ and relabeling $W_0$ as $W$, we may assume that $g$ is $C^1$ on $W$.} Next, applying \cite[Lemma~3.1]{PeterzilStarchenko2007} to each $\widehat F_i$, we define the set
\begin{equation} \label{eq:def-wstar}
W_*:= \bigcap_{i=1,\dots,m}\bigcap_{j=1,\dots,m-1} W_{ij}, \qquad 
W_{ij} := \left\{ a\in W: {\lim}_{\substack{(y,t)\to(a,0)\\(y,t)\in \mathcal U}}~\frac{\partial \widehat F_i}{\partial y_j}(y,t) = \frac{\partial g_i}{\partial y_j}(a) \right\}.
\end{equation}
The set $W_*$ is definable. By \cite[Lemma~3.1]{PeterzilStarchenko2007}, each scalar limit identity in \eqref{eq:def-wstar} holds outside a definable subset of $W$ of dimension less than $m-1$.\footnote{
{More specifically, the tangential derivative-limit lemma \cite[Lemma~3.1]{PeterzilStarchenko2007} gives, for each pair $(i,j)$,
a definable exceptional set $E_{ij}\subset W$ with
$\dim(E_{ij})<m-1$, outside which the limit identity holds.
The finite union of these exceptional sets still has
dimension less than $m-1$.}} Taking the finite union of these exceptional sets yields $\dim(W\setminus W_*)<m-1$.
Since $W$ is nonempty and open in $\mathbb{R}^{m-1}$, it follows $W_* \neq \varnothing$ and we can choose $a\in W_*$. The limits through $\mathcal U$ in \eqref{eq:def-wstar} also hold through either $\mathcal U_+$ or $\mathcal U_-$, and therefore, we obtain
\[
{\lim}_{\substack{(y,t)\to(a,0)\\(y,t)\in \mathcal U_\pm}}~D\widehat F(y,t)(\eta,0) = Dg(a)\eta \qquad \forall~\eta\in\mathbb{R}^{m-1}.
\]
Let us set $\bar z:=\Phi(a,0)$. Since $F$ is Lipschitz, the matrices $DF(z)$ is bounded whenever it is defined. The chain rule $D\widehat F(y,t) = DF(\Phi(y,t))D\Phi(y,t)$ and the continuity of $D\Phi$ then imply $\lim_{{z\to\bar z, z\in \mathcal R_\pm}} DF(z)D\Phi(a,0)(\eta,0) = Dg(a)\eta$. 
Since $\mathcal T_{\bar z}\mathcal M =D\Phi(a,0)(\mathbb{R}^{m-1}\times\{0\})$, we conclude that
\begin{equation}
 \lim_{{z\to\bar z, z\in \mathcal R_+}}DF(z)h
 =
 \lim_{{z\to\bar z, z\in \mathcal R_-}}DF(z)h
 \qquad \forall~h\in \mathcal T_{\bar z}\mathcal M.
 \label{eq:adjacent-tangential-derivative}
\end{equation}
%
%
%
Choose sequences $z_k^\pm\in \mathcal R_\pm$ converging to $\bar z$. Using the boundedness of the matrices $DF(z_k^\pm)$, $k \in \N$, and after passing to subsequences (if necessary), suppose that
\[
 A_+:=\lim_{k\to\infty}DF(z_k^+),
 \qquad
 A_-:=\lim_{k\to\infty}DF(z_k^-).
\]
The matrices $A_+$ and $A_-$ are symmetric, and nonsingular by \eqref{eq:uniform-conorm}.  In particular, the negative inertia $\nminus(A_\pm)$ is exactly the constant negative inertia on the sets $\mathcal R_\pm$.  The constancy of $z \mapsto \sgn(\det(DF(z)))$ on $\mathcal B_{F,\epsilon}(x)$, the invertibility of $A_\pm$, and the continuity of the determinant therefore yield
\begin{equation}
 \sgn(\det(A_+))=\sgn(\det(A_-)).
 \label{eq:adjacent-determinant-sign}
\end{equation}
Furthermore, by \eqref{eq:adjacent-tangential-derivative}, we have $(A_+-A_-) h= 0$ for all $h\in \mathcal T_{\bar z}\mathcal M$ and due to $\dim(\mathcal T_{\bar z}\mathcal M) = m-1$, this implies $\mathrm{rank}(A_+-A_-)\leq1$. A rank-one perturbation changes the negative inertia by at most one \cite[Corollary~4.3.9]{HornJohnson2013}, so $|\nminus(A_+)-\nminus(A_-)| \leq 1$. Since the matrices $A_+$ and $A_-$ are symmetric, and nonsingular, it holds that 
\[
 \sgn(\det(A_\pm))=(-1)^{\nminus(A_\pm)}.
\]
Equation \eqref{eq:adjacent-determinant-sign} forces $\nminus(A_+)-\nminus(A_-)$ to be even, and thus we can infer $\nminus(A_+)=\nminus(A_-)$.  Finally, connectivity of the graph $\Gamma$ implies that the negative inertia $\nminus(DF(\cdot))$ is constant on $\mathcal B_{F,\epsilon}(x)$.
\end{proof}

\begin{remark} The proof of \cref{lem:uniform-mr-gradient-inertia} is strongly inspired by the strategies used by Peterzil and Starchenko \cite[Theorem~3.2]{PeterzilStarchenko2007} and Dinh and Pham \cite[Theorem~3.1]{DinhPham}. Here, compared to \cite{PeterzilStarchenko2007,DinhPham}, we establish the constancy of the full negative inertia, rather than only the sign of the determinant. This relies on three additional ingredients: the symmetry of $DF$ due to the underlying gradient structure, the boundedness of $DF$ arising from Lipschitz continuity, and the uniform spectral gap supplied by the metric regularity of $F$.
\end{remark}

\subsection{A lineality principle for two-homogeneous mappings}

The proof combines the preceding negative-inertia lemma with a second-order
reduction to an exposed face and induction on the intrinsic dimension.

\begin{theorem}[General two-homogeneous lineality principle]
\label{thm:general-model-lineality}
Let $q:\R^n\to\Rex$ be lsc, proper, definable, $2\rho$-weakly convex, $\rho \geq 0$, and positively $2$-homogeneous.  
If $M := \partial (\frac12 q)$ is metrically regular at $(0,0)$ and
\begin{equation}
q(w)\geq0\qquad \forall~w \in\mathscr L_q := \dom(q) \cap (-\dom(q)),
\label{eq:general-model-lineality-sign}
\end{equation}
then there is $\gamma>0$ such that
\begin{equation}
q(w)\geq\gamma\norm{w}^2\qquad \forall~w\in \dom(q),
\label{eq:general-model-uniform-positive}
\end{equation}
and $M$ is strongly metrically regular at $(0,0)$.
\end{theorem}

\begin{proof}
Applying \Cref{lem:intrinsic-red}, we may assume (without loss of generality) that $\eq := \Span(\dom(q)) = \R^n$. (Indeed, as in \Cref{lem:intrinsic-red}, let the columns of $P \in \R^{n \times r}$ form an orthonormal basis of $\eq$ and set $q_\eq := q \circ P$. Noticing $\linq, \dom(q) \subset \eq$ and $\dom(q_\eq) = P^\top\dom(q)$, we may equivalently work with $q$ or $q_\eq$ in \eqref{eq:general-model-lineality-sign} and \eqref{eq:general-model-uniform-positive}.) Our proof is based on an induction on the dimension $n = \dim(\eq)$. In the following, we set $S := \dom(q)$.

Once \eqref{eq:general-model-uniform-positive} is established, the origin is a strong local minimizer of $\frac12 q$, and strong metric regularity then follows from \Cref{cor:second-order-stationarity}---applied to $h = \frac12 q$ and using $\mathrm d^{2}(\frac12 q)(0|0) = q$. The function $\frac12 q$ is weakly convex, hence prox-regular and subdifferentially continuous (at $0$ for $0$) \cite[Example~13.30 and Exercise~13.35]{RockafellarWets}. Choose $\lambda>0$ with $\lambda\rho<1$ and set
\[
 e:=\env{q/2}{\lambda},
 \qquad
 G:=\nabla e.
\]
By \Cref{prop:prox-env} and \Cref{lem:moreau-transform}~(iii), $e$ is definable, positively $2$-homogeneous, and $C^{1,1}$ on $\R^n$ and $G$ is definable, positively homogeneous, and globally Lipschitz. (Homogeneity can be verified using the definition of the Moreau envelope or as in \eqref{eq:semi-grad-env-model}). 
Moreover, \Cref{thm:exact-reduction} is applicable (with $h = \frac12 q$) and thus, $G$ is metrically regular at $(0,0)$. The mapping $M$ is metrically regular at $(0,0)$ per assumption. Consequently, there are  $\epsilon>0$ and $\kappa>0$ such that
\begin{equation}
\begin{aligned}
 \dist(z,G^{-1}(y))
 &\leq\kappa\norm{G(z)-y}
 && \forall~z\in\bB_\epsilon,\ y\in\bB_\epsilon,\\
 \dist(x,M^{-1}(y))
 &\leq\kappa\dist(y,M(x))
 && \forall~x\in\bB_\epsilon,\ y\in\bB_\epsilon.
\end{aligned}
\label{eq:general-local-mr}
\end{equation}
Next, we globalize these estimates. Since $q$ is positively $2$-homogeneous, we have $M(tx)=tM(x)$, $t>0$. Indeed, this identity follows first for regular subgradients directly from their defining inequality and then for limiting subgradients by passage to convergent graph sequences. Hence, we obtain
\[
 M^{-1}(ty)
 =\{tx:ty\in M(tx)\}
 =\{tx:y\in M(x)\}
 =tM^{-1}(y) \qquad t > 0.
\]
Fix $x,y\in\R^n$ and choose $0<t<\min\{ \frac{\epsilon}{1+\norm{x}},\frac{\epsilon}{1+\norm{y}}\}$. Then $tx,ty\in\bB_\epsilon$, and the second estimate in \eqref{eq:general-local-mr} gives
\begin{align*}
 t\dist(x,M^{-1}(y))
 =\dist(tx,M^{-1}(ty))
 \leq\kappa\dist(ty,M(tx))
 =t\kappa\dist(y,M(x)).
\end{align*}
Dividing by $t>0$ globalizes the estimate for $M$. The proof for $G$ is identical, and hence, we can infer:
\begin{equation}
 \dist(z,G^{-1}(y)) \leq\kappa\norm{G(z)-y},\qquad
 \dist(x,M^{-1}(y)) \leq\kappa\dist(y,M(x))
 \qquad \forall~x,y,z\in\R^n.
\label{eq:general-global-mr}
\end{equation}

The induction claim is immediate if $n=0$. If $n=1$, then $S = \dom(q)$ either coincides with $\R$, $\R_+$, or $\R_-$, since $0\in S$ and $S$ is a convex cone. For $S=\R$ and leveraging the positive $2$-homogeneity of $q$, we may write $q(x)=a_+x^2$ on $\R_+$ and $q(x)=a_-x^2$ on $\R_-$.
The sign assumption \eqref{eq:general-model-lineality-sign} yields $a_\pm\geq0$, and metric regularity excludes $a_+=0$ or $a_-=0$. Otherwise the range of $M$ would be contained in either $\R_-$ or $\R_+$ contradicting the local openness of $M$ around $(0,0)$ (which is ensured by metric regularity; see \Cref{sec:intro-metric}). For $S=\R_+$, write
$q(x)=ax^2+\iota_{\R_+}(x)$.  If $a\leq0$, the local range of $M$ is
contained in $\R_-$, which again contradicts local openness, and hence we have $a>0$. The case $S=\R_-$ is identical. Thus, the claim holds for dimensions $n\leq1$.

Consider the case $n\geq2$. We use the following key observation: if $G$ is $C^1$ near some $z_*\neq0$ and $DG(z_*)\succ0$, then $e$ is strongly convex. Indeed, if $G$ is differentiable at $z_*$ and using the positive homogeneity of $G$, we obtain $G(tz_* +h)-G(tz_*) = t[G(z_*+h/t)-G(z_*)]=DG(z_*)h + o(\|h\|)$ for all $h \in \Rn$ and $t > 0$. Hence, $G$ is differentiable along the ray $t \mapsto tz_*$, $t > 0$, with $DG(tz_*) = DG(z_*) \succ0$. Let $\mathcal B_{G,\epsilon^\prime}(0)$ be given as in \Cref{lem:uniform-mr-gradient-inertia} applied to $G = \nabla e$ at $(0,0)$ for some $\epsilon^\prime > 0$. By the constancy of the negative inertia in $\mathcal B_{G,\epsilon^\prime}(0)$ and the positive homogeneity of $G$, we then have $DG(z)\succeq 0$ for all $z$ where $G$ is $C^1$. Combining these results with \eqref{eq:uniform-conorm} (and \eqref{eq:general-global-mr}), we can infer $\sigma_{\min}(DG(z)) = \lambda_{\min}(\nabla^2 e(z)) \geq \kappa^{-1}$ for all such $z$. By \cite[Lemmas~2.6~(i) and~2.7]{DinhPham}, the definable mapping $G$ is $C^1$ almost everywhere and hence, applying the Hessian criterion \cite[p.~170]{RockafellarHiddenConvexity} proves the assertion.

We now verify existence of such $z_*$ depending on $S$. If $S=\R^n$, then $\frac12q\geq0$ by \eqref{eq:general-model-lineality-sign}, and hence, we have $e\geq0$. The map $G$ is open by \eqref{eq:general-global-mr} and has finite fibers by \cite[Theorem~3.1]{DinhPham}. Moreover, the condition $e(z)=0$ implies $G(z)=0$ by first-order optimality, and positive homogeneity and finite fibers force $z=0$. Consequently, it follows
\[
 c_0:={\inf}_{\norm{u}=1}~e(u)>0,
\]
and positive $2$-homogeneity gives $e(z)\geq c_0\norm{z}^2$ for all $z\in\R^n$. Let $\Sigma$ be the set of points $z$ where $G$ is not $C^1$ on any neighborhood of $z$. By \cite[Lemma~2.7 and Lemma~3.3~(i)]{DinhPham}, we have $\dim(G(\Sigma))=\dim(\Sigma)<n$. Choose $b\notin G(\Sigma)\cup\{0\}$.  At a global minimizer $z_b$ of the coercive function $e-\iprod{b}{\cdot}$, we have $G(z_b)=b$. Hence, $G$ is $C^1$
near $z_b\neq0$, and $DG(z_b) = \nabla^2 e(z_b)$ is positive semidefinite by the second-order necessary optimality conditions. Metric regularity \eqref{eq:general-global-mr} further implies that $DG(z_b)$ is nonsingular and we obtain $DG(z_b)\succ0$. Thus, $e$ is strongly convex.

Next, we consider the case where $S$ is a proper cone with $S\neq\R^n$. In the following, we introduce a reduction to expose, by a nonzero polar direction, the two-sided face of $S$ and thereby pass to a second-order model of strictly smaller intrinsic dimension. Since $q$ is weakly convex and positively $2$-homogeneous, the effective domain $S = \dom(q)$ is a convex cone. Set $C:=\cl(S)$. Moreover, $C\neq\R^n$, since $\ri(S) = \ri(C)$ \cite[Proposition~2.40]{RockafellarWets}. Hence, $C$ is a
proper closed convex cone. By the bipolar theorem \cite[Theorem~14.1]{RockafellarConvexAnalysis}, we know that $C^\circ\neq\{0\}$. Choose
$0\neq v\in\ri(C^\circ)$. Then by \cite[Theorem~12]{FenchelConvexCones}, we have 
\begin{equation*}
 C\cap \{v\}^\perp=\lin(C).
\end{equation*}
Set $S_v:=S\cap \{v\}^\perp=S\cap\lin(C)$ and $f := \frac12 q-\iprod{v}{\cdot}$. Then $\dim(\Span(S_v))\leq \dim(\lin(C))<n$ and
\begin{equation}
 S_v\cap(-S_v)=S\cap(-S)=\linq,
\label{eq:general-exposed-lineality}
\end{equation}
because $S\cap (-S)\subseteq \lin(C)$. The convex and positive $2$-homogeneous function $\frac12q+\frac\rho2\norm{\cdot}^2$ is nonnegative, and we have $\iprod{v}{x}\leq0$ for all $x\in S$. This yields
$\frac12q(x)+\frac\rho2\norm{x}^2\geq\iprod{v}{x}$ for all $x\in\R^n$, and implies $v\in\partial(\frac12 q+\frac\rho2\norm{\cdot}^2)(0)=\partial (\frac12q)(0)$ and $0\in\partial f(0)$. Positive two-homogeneity of $q$ further gives
\begin{equation}
 \Delta_t^2f(0|0)(w')
 =q(w')-\frac2t\iprod{v}{w'}.
\label{eq:general-exposed-difference-quotient}
\end{equation}
In \eqref{eq:general-exposed-difference-quotient}, we have $\langle v,w' \rangle\leq 0$ whenever $w'\in S$, and $q(w')=\infty$ if $w'\notin S$. Hence, since $q$ is lsc, taking the limit inferior as $t\downarrow0$ and $w'\to w$, we obtain
\begin{equation}
 \widehat q(w) := \rd^2f(0|0)(w)
 =q(w)+\iota_{\{v\}^\perp}(w)
 =q(w)+\iota_{S\cap\lin(C)}(w)
 =q(w)+\iota_{S_v}(w)
 \qquad \forall~w\in\R^n.
\label{eq:general-exposed-second-subderivative}
\end{equation}
The second estimate in \eqref{eq:general-global-mr} makes $\partial f=M-v$ metrically regular at $(0,0)$. Applying \Cref{thm:exact-reduction} to $f$, the model $\widehat M:=\partial(\widehat q/2)$ is metrically regular (at $(0,0)$). Using the representation $\widehat q=q+\iota_{\{v\}^\perp}$, $\widehat q$ is proper, lsc, definable, and positively two-homogeneous, and $\widehat q+\rho\norm{\cdot}^2 = (q+\rho\norm{\cdot}^2)+\iota_{\{v\}^\perp}$ is convex. Moreover, the dimension of $\dom(\widehat q)=S_v$ is strictly less than $n$. Combining \eqref{eq:general-exposed-lineality} and \eqref{eq:general-exposed-second-subderivative}, condition \eqref{eq:general-model-lineality-sign} also holds with $q$ being substituted by $\widehat q$. Thus, we can directly apply the induction hypothesis to $\widehat q$. As a result, $\widehat M$ is strongly metrically regular and there is $\gamma_v>0$ such that
\begin{equation}
 \widehat q(w)=q(w)\geq\gamma_v\norm{w}^2\qquad \forall~w\in S_v.
\label{eq:general-exposed-positive}
\end{equation}
Set $e_f:=\env{f}{\lambda}$. A direct calculation then gives
\begin{equation}
 e_f(w)=e(w+\lambda v)-\iprod{v}{w}-\frac\lambda2\norm{v}^2,
 \qquad
 \nabla e_f(w)=G(w+\lambda v)-v.
\label{eq:general-tilted-envelope-identity}
\end{equation}
In particular, $e_f$ is again definable and $C^{1,1}$ on $\R^n$. Moreover, by \Cref{thm:exact-reduction}, both $\partial f$ and $\nabla e_f$ are strongly metrically regular at $(0,0)$. Equations \eqref{eq:general-exposed-second-subderivative} and \eqref{eq:general-exposed-positive}, together with \cite[Theorem~13.24~(c)]{RockafellarWets}, make the origin a local minimizer of $f$. Due to $\prox_{\lambda f}(0)=0$ (which follows from $0 \in \partial f(0)$) and using the continuity of $\prox_{\lambda f}$, the origin is also a local minimizer of the function $e_f$. Thus, by \cite[Theorem~3.3]{DrusvyatskiyLewis}, there are $\kappa_f>0$ and neighborhoods $U_0,V$ of $0$ 
such that the following uniform quadratic-growth estimate holds
\begin{equation}
 e_f(z)\geq e_f(w)+\iprod{\nabla e_f(w)}{z-w}+\kappa_f\norm{z-w}^2,
 \quad \forall~z,w\in U_0 \,\text{with}\, \nabla e_f(w)\in V.
\label{eq:general-tilt-uniform-growth}
\end{equation}
By the continuity of $\nabla e_f$ and due to $\nabla e_f(0)=0$, we may shrink $U_0$ to a neighborhood $U$ of $0$ such that $\nabla e_f(U)\subset V$. Hence, without loss of generality, we may assume that \eqref{eq:general-tilt-uniform-growth} holds for all $z,w \in U$.
Next, choose $w_*\in U\setminus\{-\lambda v\}$ sufficiently close to zero such that $\nabla e_f$ is $C^1$ on a neighborhood of $w_*$; this is possible by \cite[Lemma~2.7]{DinhPham} since $\nabla e_f$ is definable. Set $z_*:=w_*+\lambda v$. By \eqref{eq:general-tilted-envelope-identity}, $G$ is $C^1$ on a neighborhood of $z_*$. Applying \eqref{eq:general-tilt-uniform-growth} at $w_*$ and using a second-order Taylor expansion then yields
\[
 DG(z_*)=\nabla^2e_f(w_*)\succeq2\kappa_f I\succ0.
\]
Therefore, we can again infer that $e$ is strongly convex on $\R^n$.

In either case, strong convexity, $e(0)=0$, and $G(0)=0$ give $e(x)\geq c\norm{x}^2$ for some $c>0$.  Since $e(x)\leq \frac12q(x)$ for all $x\in \R^n$, this proves \eqref{eq:general-model-uniform-positive} and finishes the induction. Thus, the estimate \eqref{eq:general-model-uniform-positive} holds in every dimension and the origin is a local minimizer of $\frac12q$. By \cite[Theorem~13.24~(a)]{RockafellarWets}, we then have $\rd^2(\frac12 q)(0|0)(w)\geq0$ for all $w\in\R^n$. Consequently, \Cref{cor:second-order-stationarity} shows that $M=\partial (\frac12 q)$ is strongly metrically regular at $(0,0)$.
\end{proof}

\Cref{thm:exact-reduction} allows us to transfer strong metric regularity of $\partial(\frac12 q)$ back to the original subgradient mapping $\partial h$. 

\begin{theorem}[A second-subderivative-based lineality criterion for strong metric regularity]
\label{thm:general-subgradient-lineality}
Let $h:\R^n\to\Rex$ be proper, lsc, and definable, and let $(\bar x,\bar v)\in\gph(\partial h)$ be given. Assume that $h$ is prox-regular and subdifferentially continuous at $\bar x$ for $\bar v$, and that
$\partial h$ is metrically regular at $(\bar x,\bar v)$. If
\[ q(w) := \rd^2h(\bar x|\bar v)(w) \geq 0 \qquad \forall~w \in \mathscr L_q:=\dom(q)\cap(-\dom(q)).
\]
then there is $\gamma > 0$ such that $q(w) \geq \gamma \|w\|^2$ for all $w \in \dom(q)$ and $\partial h$ is strongly metrically regular at $(\bar x,\bar v)$.
\end{theorem}

\begin{proof}
By \Cref{prop:automatic-twice-epi}, $q$ is proper. As discussed in \Cref{sec:basicprox}, $q$ is lsc, positively $2$-homogeneous, and $2\rho$-weakly convex (by \cite[Propositions~13.5 and 13.49]{RockafellarWets}). By \Cref{thm:exact-reduction}, $M=\partial(\half q)$ is metrically regular at $(0,0)$. \Cref{thm:general-model-lineality} implies that $q$ satisfies the asserted quadratic growth property \eqref{eq:general-model-uniform-positive} and $M$ is strongly metrically regular. Invoking \Cref{thm:exact-reduction} once more, this yields strong metric regularity of $\partial h$ (at $(\bar x,\bar v)$).
\end{proof}

As an immediate application of \Cref{thm:general-subgradient-lineality}, we obtain the following result. 

\begin{corollary}[Pointed domains]
\label{cor:general-pointed-domain}
Let $h : \Rn \to \Rex$ be proper, lsc, definable, and prox-regular and subdifferentially continuous at some $\bar x$ for $\bar v \in \partial h(\bar x)$. Setting $q := \rd^2h(\bar x|\bar v)$, suppose that the effective domain $\dom(q)$ is pointed, i.e., $\linq := \dom(q) \cap (-\dom(q)) = \{0\}$. Then, metric regularity of $\partial h$ at
$(\bar x,\bar v)$ implies strong metric regularity there.
\end{corollary}


\section{Conclusion and future directions}
\label{sec:conclusion}

We have given two complementary answers to when metric regularity of a
definable prox-regular subdifferential $\partial h$ upgrades to strong metric regularity. A Moreau envelope-based reduction transfers this question to the positively homogeneous gradient model $H = \nabla \env{q/2}{\lambda}$ determined by the second subderivative $q = \rd^2 h(\bar x|\bar v)$.
For this model, the degree formula $\deg_0(H)=1-\chi(\negq)$ provides a topological link to the set of (normalized) negative curvature directions of $q$ and yields a necessary and sufficient criterion to exclude branching. Furthermore, we present a lineality principle that gives a geometric sufficient condition: nonnegativity of the second subderivative $q$ on the lineality space $\dom(q)\cap(-\dom(q))$ propagates to uniform quadratic positivity on the whole effective domain. In particular, metric regularity implies strong metric regularity at (weak and standard) second-order stationary points.



\paragraph{Weak saddle point avoidance.} The lineality principle established in \Cref{sec:lineality} motivates the following concepts, properties, and broader research questions. 

\begin{definition}[Strong saddle points and weak saddle point avoidance] \label[definition]{def:weak-saddle} Let $h : \Rn \to \Rex$ be a proper, lsc function and let $\bar x$ be a stationary point of $h$. Set $q := \rd^2 h(\bar x|0)$. We call $\bar x$ a \emph{strong saddle point} of $h$ if
\[ \exists~w \in \dom(q) \cap (-\dom(q)) \quad \text{such that} \quad q(w) = \rd^2 h(\bar x|0)(w) < 0. \]
We say that $h$ satisfies the \emph{weak saddle property} if every stationary point is either a local minimizer or a strong saddle point. 
An algorithm $\mathscr A$ for solving $\min_x h(x)$ satisfies the \emph{weak saddle point avoidance} property (under certain algorithmic assumptions) if every accumulation point of its iterates is a weak second-order stationary point of $h$. 
\end{definition}

Under definability, prox-regularity, and subdifferential continuity, \Cref{thm:general-subgradient-lineality} shows that every weak second-order stationary point at which $\partial h$ is metrically regular is a strong local minimizer and that $\partial h$ is strongly metrically regular there. Consequently, if metric regularity holds at every stationary point, then $h$ satisfies the weak saddle property. Weak saddle point avoidance would therefore yield local minimality and strong metric regularity at all accumulation points. More generally, it is natural to ask which structural conditions on $h$ ensure the weak saddle property without requiring metric regularity at every stationary point. Several additional algorithmic questions remain open: 
\begin{itemize}
\item Can the localized Moreau envelope-based reduction and the positively homogeneous second-order models developed here be used to establish instability at strong saddle points? 
\item Which additional assumptions ensure that the attraction regions of strong saddle points are negligible, or that suitable stochastic perturbations prevent convergence to them? 
\item Can proximal, (stochastic) subgradient, or curvature-based methods satisfy weak saddle point avoidance when the second subderivative $q$ is not generalized quadratic?
\end{itemize}
These questions are closely related to the active strict saddle framework of Davis and Drusvyatskiy \cite{DavisDrusvyatskiy}. An active strict saddle admits an active $C^2$-manifold on which the objective has a negative tangential curvature direction (together with a sharpness condition). This structure enables a smooth local analysis of proximal updates and leads to saddle avoidance results for proximal methods \cite{DavisDrusvyatskiy} and, under additional assumptions, stochastic subgradient methods \cite{BianchiHachemSchechtman}. In the $C^2$-active manifold setting, the second-order model $q$ is generalized quadratic, so active strict saddles are strong saddle points in the sense of \Cref{def:weak-saddle}; see also \cite[Lemma 4.12]{qiu2026saddle}. It is interesting to investigate whether and in which structural and algorithmic scenarios avoidance beyond this smooth active manifold setting is possible---allowing more general non-quadratic second-order models. Whether the geometric and topological results presented here can support such an analysis remains an intriguing question. We leave these directions for future work.

\section*{Acknowledgments}

During the preparation of this work, the authors used OpenAI Codex, a
generative-AI system, as an assistant for mathematical reasoning, including
proof checking and development, and for manuscript preparation. The authors have independently reviewed, verified, and essentially rewritten all mathematical claims and take full responsibility for the content of the paper. Below, we discuss some of our key ideas, approaches, and steps in more detail.

We initially tried to establish the implication in \eqref{eq:main-problem} directly at general reference pairs by combining definable stratifications with Riemannian Hessians on the smooth strata. (This idea is motivated by the deep connections between subgradients and smooth Riemannian gradients \cite{bolte2007clarke}). However, this approach was not very fruitful; a major difficulty is to control the second-order behavior across stratum boundaries. At this point, Codex also proposed a semi-algebraic counter-example---illustrating that metric and strong metric regularity of the subdifferential $\partial h$ are not automatically and generally equivalent in the definable setting. (Here, let us note that the earlier preprint versions of \cite{MaSun} do not contain such a counter-example; it was added in the most recent version~5, posted on 10~August~2026 \cite{MaSun}).

Motivated by \cite{DMN,chen2025aubin,chen2026characterization}, we then sought proper characterizations of \eqref{eq:main-problem} at more specific reference points, such as local minimizers. Based on the variational second-order perspective adopted in our earlier work \cite{OuyangMilzarekI,OuyangMilzarekII}, we first realized that definability and prox-regularity yield proper twice epi-differentiability, as shown in \Cref{prop:automatic-twice-epi}. This key (but somewhat folkloric) observation led to the various reduction steps presented in \Cref{sec:model-reduction}---through a localized Moreau envelope of $h$, the second subderivative $q=\rd^2h(\bar x|\bar v)$ of $h$, and the Moreau envelope of $q/2$. The weak convexity and positive two-homogeneity of $q$ yield a globally Lipschitz continuous, positively homogeneous gradient $\nabla\env{q/2}{\lambda}$ for sufficiently small $\lambda>0$, which enables a degree-theoretic analysis under metric regularity (not just for local minimizer but for general reference points). This route provided a more powerful and rich link to variational analysis compared to our previous stratification-based attempts. Codex helped to formulate some of the auxiliary lemmas and to develop, check, and revise the arguments relating this degree to the Euler characteristic of $\negq$.

As the Euler characteristic-based curvature criterion (\Cref{thm:abstract-euler}) already fully addresses the local-minimizer setting, we were wondering whether the results can be pushed further to weak second-order stationary points and if non-negativity of $q$ on the lineality space of $\dom(q)$ is sufficient to ensure the implication in \eqref{eq:main-problem}. Codex helped to develop and test some of the proof arguments for the resulting lineality principle, including the negative inertia propagation and the induction step used in the proof of \Cref{thm:general-model-lineality}. Throughout this process, we corrected candidate arguments, supplied (many) further proof steps and details, and (essentially completely) restructured the presentation and discussions.


\appendix 

\section{Completing the proof of \Cref{lem:homogeneous-index}} \label[appendix]{sec:app:remainder-of-lemma}

\begin{proof}
In order to complete the proof of \Cref{lem:homogeneous-index}, we need to show $\chi(\{\phi<0\})=\chi(\{\widetilde\phi<0\})$.  
First, we note that the strict estimates \eqref{eq:homotopy-field-lower-bound} and \eqref{eq:zero-level-gradient-bound} remain valid for all $u\in\sphL$ and $s$ in some open interval $I\supseteq[0,1]$. Set
\begin{align*}
\mathcal X&:=\{(s,u)\in I\times\sphL:\phi_s(u)\le0\}, \quad \mathcal X_-:=\{(s,u)\in\mathcal X:\phi_s(u)<0\},\quad
\mathcal X_0:=\{(s,u)\in\mathcal X:\phi_s(u)=0\},
\end{align*}
and $F:I\times\sphL\to\R$, $F(s,u):=\phi_s(u)$. Our strategy is to show that $\mathcal X_-$ and $\mathcal X_0$ form a $C^1$-Whitney stratification of $\mathcal X$. This allows us to apply Shiota's definable first isotopy lemma and to verify that $\{\phi<0\}$ and $\{\widetilde\phi<0\}$ are homeomorphic.

The map $F$ is definable and $C^1$, so $\mathcal X_-=F^{-1}((-\infty,0))$ is a definable open $C^1$-submanifold of $I\times\sphL$. Moreover, using \eqref{eq:zero-level-gradient-bound}, we have
\[
\rd_uF(s,u)=\rd\phi_s(u)\ne0 \quad \forall~(s,u)\in\mathcal X_0
\]
Thus, $0$ is a regular value of $F$, and the $C^1$-regular-value theorem \cite[Chapter~1, Theorems~3.2--3.3]{Hirsch} shows that $\mathcal X_0=F^{-1}(0)$ is a definable $C^1$-submanifold of $I\times\sphL$. Consequently, $\{\mathcal X_-,\mathcal X_0\}$ is a finite definable $C^1$-stratification of $\mathcal X$ in the sense of \cite[Definition~2.3]{Tanabe2023}.
The frontier condition \cite[Definition~2.9]{Tanabe2023} holds as well.
Indeed, by \cref{eq:zero-level-gradient-bound} (and its extension to $s\in I$), we obtain for all $(s,u)\in\mathcal X_0$
\[
\rd\phi_s(u)[v]
=-\norm{\nabla_{\sphL}\phi_s(u)}^2<0,
\]
where $v:= -\nabla_{\sphL}\phi_s(u)\in \mathcal T_u\sphL$. Choose a $C^1$-curve $\gamma:(-\varepsilon,\varepsilon)\to\sphL$ satisfying
$\gamma(0)=u$ and $\gamma'(0)=v$. Then, we have $\phi_s(\gamma(t)) = \phi_s(u) + t \rd\phi_s(u)[v] + o(t) = -t \norm{\nabla_{\sphL}\phi_s(u)}^2 + o(t) <0$ for all sufficiently small $t>0$ and hence, it follows $\mathcal X_- \ni (s,\gamma(t)) \to (s,u)$ as $t \downarrow 0$. This shows $\mathcal X_0 \subset \cl({\mathcal X_-})$ and we have
\[
\cl({\mathcal X_-})\cap\mathcal X_0=\mathcal X_0,
\qquad
\cl({\mathcal X_0}) \cap\mathcal X_-=\varnothing.
\]
It remains to verify Whitney $(b)$-regularity. Using its sequential definition \cite[Section~2, Definition~2.1]{Mather2012}, let $x_k=(s_k,u_k)\in\mathcal X_-$ and $y_k=(\widehat s_k,\widehat u_k)\in\mathcal X_0$ converge to $y=(s,u)\in\mathcal X_0$. Suppose that $\mathcal T_{x_k}\mathcal X_-$ converges to a linear subspace $\mathcal T$ and that the secant lines $L_k:=\R(x_k-y_k)$ converge to a line $\ell$. Since $\mathcal X_-$ is open in $I\times\sphL$, we have 
\[
\mathcal T_{x_k}\mathcal X_-
=\R\times \mathcal T_{u_k}\sphL
=\R\times\{z\in\sL:\iprod{z}{u_k}=0\}
\rightarrow \R\times \mathcal T_u\sphL,
\]
which implies $\mathcal T =\R\times \mathcal T_u\sphL$.  
After passing to a subsequence, the unit vectors ${(x_k-y_k)}/{\|x_k-y_k\|}\in L_k$ converge to a unit vector $(a,z)$ spanning $\ell$.
Since $\norm{u_k}=\norm{\widehat u_k}=1$, we obtain
\[ \frac{\iprod{u_k+\widehat u_k}{{u_k-\widehat u_k}}}{{\|x_k-y_k\|}} = 0 \]
Passing to the limit shows that any limiting unit secant $(a,z)$ satisfies $\iprod{u}{z}=0$. Therefore, $\ell \subset \mathcal T = \lim_{k \to \infty} \mathcal T_{x_k}\mathcal X_-$. 
This verifies Whitney $(b)$-regularity. Hence, by \cite[Definition~2.10]{Tanabe2023}, the nonempty members of $\{\mathcal X_-,\mathcal X_0\}$ form a finite definable $C^1$-Whitney stratification of $\mathcal X$.

Next, we verify that the parameter projection 
\[
\pi:\mathcal X\to I, \qquad \pi(s,u):=s,
\]
is proper (i.e., compact parameter sets have compact inverse images under $\pi$) and $\pi$ is a submersion on each stratum (i.e., its differential is surjective). These properties are the remaining hypotheses required by the first isotopy lemma.
%
%
%
%
The projection $\pi$ is clearly proper, since $\pi^{-1}(K)=\mathcal X\cap (K\times \sphL)=\{(s,u)\in K\times \sphL:~\phi_s(u)\leq 0\}$ is a closed subset of the compact set $K\times \sphL$ for any compact set $K\subseteq I$.  It is also a submersion on both strata: on $\mathcal X_0$, this follows from $\rd\phi_s(u)\ne0$
\footnote{At $(s,u)\in\mathcal X_0$, we have $\mathcal T_{(s,u)}\mathcal X_0 =\{(\tau,v)\in\R\times \mathcal T_u\sphL: \tau\,\partial_s\phi_s(u)+\rd\phi_s(u)[v]=0\}$.
For every $\tau\in\R$, since $\rd\phi_s(u)\ne0$, we can choose $v\in \mathcal T_u\sphL$ such that $\rd\phi_s(u)[v]=-\tau\,\partial_s\phi_s(u)$.  Hence $D(\pi|_{\mathcal X_0})(s,u)(\tau,v)=\tau$, so $D(\pi|_{\mathcal X_0})(s,u)$ is surjective.}, and on $\mathcal X_-$ this is clear from the openness. Finally, we verify that each nonempty stratum is projected onto $I$. Define
\[
\alpha(s):={\min}_{u\in\sphL}~\phi_s(u),
\qquad s\in I.
\]
The compactness of $\sphL$ and the identity $\phi_s-\phi_t=(s-t)(\widetilde\phi-\phi)$ yield
\[
|\alpha(s)-\alpha(t)|
\le \norm{\phi_s-\phi_t}_\infty
=|s-t|\norm{\widetilde\phi-\phi}_\infty.
\]
Thus, $\alpha$ is continuous and we have
\[
\pi(\mathcal X_-) = \{s \in I: \exists~u \in \sphL\,\text{such that}\,\phi_s(u) < 0\} = \{s\in I:\alpha(s)<0\}.
\]
Suppose that $\mathcal X_-$ is nonempty and that $A:=\pi(\mathcal X_-)$ is a proper subset of $I$. The equation above and the continuity of $\alpha$ show that $A$ is open in $I$. Since $I$ is an interval, it is connected \cite[Corollary~24.2]{Munkres}. Thus, the nonempty proper subset $A$ cannot be both open and closed in $I$ \cite[Section~23, p.~148]{Munkres}, and its boundary relative to $I$ is nonempty \cite[Section~17, Exercise~19(b)]{Munkres}. Choose a point $s_*\in I$ in this relative boundary. Since $A$ is open, $s_*\notin A$, and hence $\alpha(s_*)\geq0$. Furthermore, there are points $s_k\in A$ with $s_k\to s_*$.  Therefore, $\alpha(s_k)<0$ for every $k$, and the continuity of $\alpha$ gives $\alpha(s_*)\leq0$. Consequently, $\alpha(s_*)=0$. At a point $u_*\in\sphL$, where $\phi_{s_*}$ attains this minimum, we have
\[
\phi_{s_*}(u_*)=0,
\qquad
\rd\phi_{s_*}(u_*)=0.
\]
This implies $V_{\phi_{s_*}}(u_*)=0$, contradicting \eqref{eq:homotopy-field-lower-bound}. Hence, either $\mathcal X_-$ is empty or $\pi(\mathcal X_-)=I$. Moreover, $\pi(\mathcal X_0)$ is open in $I$ because $\pi|_{\mathcal X_0}$ is a submersion and hence an open map \cite[Proposition~4.1.2]{MargalefOuterelo}, and it is closed in $I$ because $\mathcal X_0$ is closed in $I\times\sphL$ and the coordinate projection $I\times\sphL\to I$ is a closed map \cite[Section~26, Exercise~7]{Munkres}. Thus, the connectedness of $I$ gives $\pi(\mathcal X_0)=I$ whenever $\mathcal X_0$ is nonempty.

If $\mathcal X_-$ is empty, then both $\{\phi<0\}$ and $\{\widetilde\phi<0\}$ are empty, and their Euler characteristics are 
zero.  Otherwise, Shiota's definable first isotopy lemma, in the form
stated in \cite[Theorem~2.14]{Tanabe2023}, shows that $\pi_-:=\pi|_{\mathcal X_-}:\mathcal X_-\to I$ is a definable $C^1$-locally trivial fibration \cite[Definition~2.13]{Tanabe2023}, i.e., for every $s\in I$ there are
an open interval $S\subset I$ containing $s$ and a definable $C^1$-mapping $\rho_{s}$ such that
\[
(\pi_-,\rho_{s}):\pi_-^{-1}(S)
\rightarrow S\times\pi_-^{-1}(s)
\]
is a definable $C^1$-diffeomorphism. Consequently, the fibers $\pi_-^{-1}(\sigma)=\{\sigma\}\times\{\phi_\sigma<0\}$ are definably diffeomorphic for any two parameters in the same trivializing interval $S$, because they are all diffeomorphic to $\pi_{-}^{-1}(s)$. The intervals $S$ ($s \in I$) form an open cover of $[0,1]$ and (thanks to compactness), we can choose a finite subcover. Select a subdivision $0=s_0 < \dots < s_N=1$ such that each pair of successive parameters belongs to a common trivializing interval. The corresponding trivializations $\rho_{s_j}$ give definable $C^1$-diffeomorphisms $\theta_j : \pi_-^{-1}(s_{j-1}) \to \pi_-^{-1}(s_{j})$. Their finite composition $\theta_N\circ\cdots\circ\theta_1$ identifies the endpoint fibers 
\[ \pi_-^{-1}(0) = \{0\} \times \{\phi < 0\} \quad \text{and} \quad \pi_-^{-1}(1) = \{1\} \times \{\widetilde\phi < 0\}. \] 
Identifying $\pi_-^{-1}(s)=\{s\}\times\{\phi_s<0\}$ with $\{\phi_s<0\}$, we conclude that $\{\phi < 0\}$ and $\{\widetilde \phi < 0\}$ are definably $C^1$-diffeomorphic and hence have the same Euler characteristic (by \cite[Corollary~2.11 and Theorem~2.44]{Hatcher}).
\end{proof}

\bibliographystyle{siam}
\bibliography{references}

@article {bolte2007clarke,
    AUTHOR = {Bolte, J\'er\^ome and Daniilidis, Aris and Lewis, Adrian and
              Shiota, Masahiro},
     TITLE = {Clarke subgradients of stratifiable functions},
   JOURNAL = {SIAM J. Optim.},
  FJOURNAL = {SIAM Journal on Optimization},
    VOLUME = {18},
      YEAR = {2007},
    NUMBER = {2},
     PAGES = {556--572},
      ISSN = {1052-6234,1095-7189},
   MRCLASS = {49J52 (26D10 32B20)},
  MRNUMBER = {2338451},
MRREVIEWER = {Marian\ Mure\c san},
       DOI = {10.1137/060670080},
       URL = {https://doi.org/10.1137/060670080},
}

@article {dries1996geometric,
    AUTHOR = {van den Dries, Lou and Miller, Chris},
     TITLE = {Geometric categories and o-minimal structures},
   JOURNAL = {Duke Math. J.},
  FJOURNAL = {Duke Mathematical Journal},
    VOLUME = {84},
      YEAR = {1996},
    NUMBER = {2},
     PAGES = {497--540},
      ISSN = {0012-7094,1547-7398},
   MRCLASS = {32B20 (03C68)},
  MRNUMBER = {1404337},
MRREVIEWER = {Anand\ Pillay},
       DOI = {10.1215/S0012-7094-96-08416-1},
       URL = {https://doi.org/10.1215/S0012-7094-96-08416-1},
}

@article {bolte2009semismooth,
    AUTHOR = {Bolte, J\'er\^ome and Daniilidis, Aris and Lewis, Adrian},
     TITLE = {Tame functions are semismooth},
   JOURNAL = {Math. Program.},
  FJOURNAL = {Mathematical Programming. A Publication of the Mathematical
              Programming Society},
    VOLUME = {117},
      YEAR = {2009},
    NUMBER = {1-2},
     PAGES = {5--19},
      ISSN = {0025-5610,1436-4646},
   MRCLASS = {49J52 (03C64 03C98 14P10 65H10 90C31)},
  MRNUMBER = {2421297},
       DOI = {10.1007/s10107-007-0166-9},
       URL = {https://doi.org/10.1007/s10107-007-0166-9},
}

@phdthesis{milzarek2016numerical,
	AUTHOR = {Milzarek, Andre},
  	TITLE = {Numerical methods and second order theory for nonsmooth problems},
  	YEAR = {2016},
  	SCHOOL = {Technische Universit{\"a}t M{\"u}nchen}
}

@unpublished{BenevieriFuriPeraSpadini,
    	AUTHOR = {Benevieri, P. and Furi, M. and Pera, M.P. and Spadini, M.},
    	TITLE = {An introduction to topological degree in Euclidean spaces},
    	NOTE = {arXiv preprint, arXiv:2304.06463},
    	YEAR = {2023},
}

@unpublished{OuyangMilzarekII,
    author = {Ouyang, W. and Milzarek, A.},
    title = {Variational properties of decomposable functions. {P}art {II}:
{S}trong second-order theory},
    note = {arXiv preprint, arXiv:2311.07276v3},
    year = {2025},
}

@unpublished{qiu2026saddle,
    author = {Qiu, J. and Ma, B. and Milzarek, A. and Zhang, J.},
    title = {Stochastic Saddle Avoidance Beyond Unit Excitation and
Smoothness: A Pathwise {L}yapunov–{P}erron Framework},
    note = {arXiv preprint, arXiv:2608.03001},
    year = {2026},
}

@unpublished{OuyangMilzarekI,
    author = {Ouyang, W. and Milzarek, A.},
    title = {Variational properties of decomposable functions. {P}art {I}:
{S}trict epi-calculus and applications},
    note = {arXiv preprint, arXiv:2311.07267v2},
    year = {2024},
}

@article {BarretoFenilleHartmann,
    AUTHOR = {Barreto, A. P. and Fenille, M. C. and Hartmann, L.},
     TITLE = {Inverse mapping theorem and local forms of continuous
              mappings},
   JOURNAL = {Topology Appl.},
  FJOURNAL = {Topology and its Applications},
    VOLUME = {197},
      YEAR = {2016},
     PAGES = {10--20},
      ISSN = {0166-8641,1879-3207},
   MRCLASS = {58C15 (46T20 47J07 55M25 57N45)},
  MRNUMBER = {3426904},
MRREVIEWER = {Daciberg\ Lima\ Gon\c calves},
       DOI = {10.1016/j.topol.2015.10.013},
       URL = {https://doi.org/10.1016/j.topol.2015.10.013},
}

@article {chen2026characterization,
    AUTHOR = {Chen, Liang and Chen, Ruoning and Sun, Defeng and Zhang,
              Liping},
     TITLE = {Characterizations of the {A}ubin property of the solution
              mapping for nonlinear semidefinite programming},
   JOURNAL = {Math. Program.},
  FJOURNAL = {Mathematical Programming},
    VOLUME = {215},
      YEAR = {2026},
    NUMBER = {1-2},
     PAGES = {637--668},
      ISSN = {0025-5610,1436-4646},
   MRCLASS = {49J53 (90C22 90C31 90C46)},
  MRNUMBER = {5010065},
       DOI = {10.1007/s10107-025-02231-2},
       URL = {https://doi.org/10.1007/s10107-025-02231-2},
}

@article {hang2024strict,
    AUTHOR = {Hang, Nguyen T. V. and Sarabi, M. Ebrahim},
     TITLE = {A chain rule for strict twice epi-differentiability and its
              applications},
   JOURNAL = {SIAM J. Optim.},
  FJOURNAL = {SIAM Journal on Optimization},
    VOLUME = {34},
      YEAR = {2024},
    NUMBER = {1},
     PAGES = {918--945},
      ISSN = {1052-6234,1095-7189},
   MRCLASS = {49J52 (49J53 90C31)},
  MRNUMBER = {4711228},
       DOI = {10.1137/22M1520025},
}

@article {poliquin1996hessian,
    AUTHOR = {Poliquin, R. A. and Rockafellar, R. T.},
     TITLE = {Generalized {H}essian properties of regularized nonsmooth
              functions},
   JOURNAL = {SIAM J. Optim.},
  FJOURNAL = {SIAM Journal on Optimization},
    VOLUME = {6},
      YEAR = {1996},
    NUMBER = {4},
     PAGES = {1121--1137},
      ISSN = {1052-6234},
   MRCLASS = {49J52 (26B05 58C20)},
  MRNUMBER = {1416532},
       DOI = {10.1137/S1052623494279316},
}

@article {robinson1980regular,
    AUTHOR = {Robinson, Stephen M.},
     TITLE = {Strongly regular generalized equations},
   JOURNAL = {Math. Oper. Res.},
  FJOURNAL = {Mathematics of Operations Research},
    VOLUME = {5},
      YEAR = {1980},
    NUMBER = {1},
     PAGES = {43--62},
      ISSN = {0364-765X,1526-5471},
   MRCLASS = {90C48 (58C15 90C31)},
  MRNUMBER = {561153},
MRREVIEWER = {M.\ B.\ Suryanarayana},
       DOI = {10.1287/moor.5.1.43},
       URL = {https://doi.org/10.1287/moor.5.1.43},
}

@article{robinson1992normal,
    AUTHOR = {Robinson, Stephen M.},
     TITLE = {Normal maps induced by linear transformations},
   JOURNAL = {Math. Oper. Res.},
  FJOURNAL = {Mathematics of Operations Research},
    VOLUME = {17},
      YEAR = {1992},
    NUMBER = {3},
     PAGES = {691--714},
       DOI = {10.1287/moor.17.3.691},
       URL = {https://doi.org/10.1287/moor.17.3.691},
}

@article{ralph1993homeomorphism,
    AUTHOR = {Ralph, Daniel},
     TITLE = {A new proof of {Robinson}'s homeomorphism theorem for
              {PL}-normal maps},
   JOURNAL = {Linear Algebra Appl.},
  FJOURNAL = {Linear Algebra and its Applications},
    VOLUME = {178},
      YEAR = {1993},
     PAGES = {249--260},
       DOI = {10.1016/0024-3795(93)90344-N},
       URL = {https://doi.org/10.1016/0024-3795(93)90344-N},
}

@article{ioffe2018polyhedral,
    AUTHOR = {Ioffe, Alexander D.},
     TITLE = {On variational inequalities over polyhedral sets},
   JOURNAL = {Math. Program.},
  FJOURNAL = {Mathematical Programming},
    VOLUME = {168},
      YEAR = {2018},
    NUMBER = {1--2},
     PAGES = {261--278},
       DOI = {10.1007/s10107-016-1077-4},
       URL = {https://doi.org/10.1007/s10107-016-1077-4},
}

@article {chen2025aubin,
    AUTHOR = {Chen, Liang and Chen, Ruoning and Sun, Defeng and Zhu,
              Junyuan},
     TITLE = {Aubin property and strong regularity are equivalent for
              nonlinear second-order cone programming},
   JOURNAL = {SIAM J. Optim.},
  FJOURNAL = {SIAM Journal on Optimization},
    VOLUME = {35},
      YEAR = {2025},
    NUMBER = {2},
     PAGES = {712--738},
      ISSN = {1052-6234,1095-7189},
   MRCLASS = {49J53 (90C30 90C31 90C46)},
  MRNUMBER = {4889303},
       DOI = {10.1137/24M1670676},
       URL = {https://doi.org/10.1137/24M1670676},
}

@article {HohLabObe20,
    	AUTHOR = {Hoheisel, Tim and Laborde, Maxime and Oberman, Adam},
     	TITLE = {A regularization interpretation of the proximal point method for weakly convex functions},
   	JOURNAL = {J. Dyn. Games},
  	FJOURNAL = {Journal of Dynamics and Games},
    	VOLUME = {7},
      	YEAR = {2020},
    	NUMBER = {1},
     	PAGES = {79--96},
      	ISSN = {2164-6066},
   	MRCLASS = {65K10 (26B25 90C25)},
  	MRNUMBER = {4155077},
	MRREVIEWER = {Letizia Pellegrini},
       	DOI = {10.3934/jdg.2020005},
       	URL = {https://doi.org/10.3934/jdg.2020005},
}

@article {gfrerer2025regular,
    AUTHOR = {Gfrerer, Helmut},
     TITLE = {On second-order variational analysis of variational convexity
              of prox-regular functions},
   JOURNAL = {Set-Valued Var. Anal.},
  FJOURNAL = {Set-Valued and Variational Analysis. Theory and Applications},
    VOLUME = {33},
      YEAR = {2025},
    NUMBER = {1},
     PAGES = {Paper No. 8, 26},
      ISSN = {1877-0533,1877-0541},
   MRCLASS = {49J53 (49J52 90C31)},
  MRNUMBER = {4867341},
MRREVIEWER = {Pedro\ P\'erez-Aros},
       DOI = {10.1007/s11228-025-00744-8},
       URL = {https://doi.org/10.1007/s11228-025-00744-8},
}

@unpublished{Coste,
    	AUTHOR = {Coste, M.},
    	TITLE = {An Introduction to O-minimal Geometry},
    	NOTE = {Dottorato di Ricerca in Matematica, Dipartimento di Matematica, Universit\`a di Pisa, Istituti Editoriali e Poligrafici Internazionali, Pisa},
    	YEAR = {2000},
}

@book {DincaMawhin,
    	AUTHOR = {Dinca, George and Mawhin, Jean},
     	TITLE = {Brouwer degree---the core of nonlinear analysis},
    	SERIES = {Progress in Nonlinear Differential Equations and their Applications},
    	VOLUME = {95},
 	PUBLISHER = {Birkh\"auser/Springer, Cham},
        YEAR = {2021},
     	PAGES = {xix+447},
      	ISBN = {978-3-030-63229-8; 978-3-030-63230-4},
   	MRCLASS = {58C30 (34C25 39A23 47H11 54H25 55M20)},
  	MRNUMBER = {4390178},
       	DOI = {10.1007/978-3-030-63230-4},
       	URL = {https://doi.org/10.1007/978-3-030-63230-4},
}

@article {DinhPham,
	    	AUTHOR = {Dinh, Si Tiep and Pham, Tien-Son},
     	TITLE = {Definable continuous mappings and {W}hyburn's conjecture},
   	JOURNAL = {Proc. Amer. Math. Soc.},
  	FJOURNAL = {Proceedings of the American Mathematical Society},
    	VOLUME = {151},
      	YEAR = {2023},
    	NUMBER = {5},
     	PAGES = {2081--2095},
      	ISSN = {0002-9939,1088-6826},
   	MRCLASS = {26B10 (03C64 54C10)},
  	MRNUMBER = {4556202},
	MRREVIEWER = {Milan\ Matejdes},
       	DOI = {10.1090/proc/16272},
        URL = {https://doi.org/10.1090/proc/16272},
}

@article {PeterzilStarchenko2007,
        AUTHOR = {Peterzil, Ya'acov and Starchenko, Sergei},
         TITLE = {Computing o-minimal topological invariants using
                  differential topology},
       JOURNAL = {Trans. Amer. Math. Soc.},
        FJOURNAL = {Transactions of the American Mathematical Society},
        VOLUME = {359},
          YEAR = {2007},
        NUMBER = {3},
         PAGES = {1375--1401},
          ISSN = {0002-9947},
           DOI = {10.1090/S0002-9947-06-04220-6},
           URL = {https://doi.org/10.1090/S0002-9947-06-04220-6},
}

@article {DontchevRockafellar,
    	AUTHOR = {Dontchev, A. L. and Rockafellar, R. T.},
     	TITLE = {Characterizations of strong regularity for variational inequalities over polyhedral convex sets},
   	JOURNAL = {SIAM J. Optim.},
  	FJOURNAL = {SIAM Journal on Optimization},
    	VOLUME = {6},
      	YEAR = {1996},
    	NUMBER = {4},
     	PAGES = {1087--1105},
      	ISSN = {1052-6234},
   	MRCLASS = {90C31 (49J40 49J52 49K40)},
  	MRNUMBER = {1416530},
	MRREVIEWER = {Adam\ B.\ Levy},
       	DOI = {10.1137/S1052623495284029},
       	URL = {https://doi.org/10.1137/S1052623495284029},
}

@article {moreau1965proximite,
    AUTHOR = {Moreau, Jean-Jacques},
     TITLE = {Proximit\'e{} et dualit\'e{} dans un espace hilbertien},
   JOURNAL = {Bull. Soc. Math. France},
  FJOURNAL = {Bulletin de la Soci\'et\'e{} Math\'ematique de France},
    VOLUME = {93},
      YEAR = {1965},
     PAGES = {273--299},
      ISSN = {0037-9484},
   MRCLASS = {46.15 (41.60)},
  MRNUMBER = {201952},
MRREVIEWER = {I.\ G.\ Amemiya},
       URL = {http://www.numdam.org/item?id=BSMF_1965__93__273_0},
}

@book {dontchev2009implicit,
    AUTHOR = {Dontchev, Asen L. and Rockafellar, R. Tyrrell},
     TITLE = {Implicit functions and solution mappings},
    SERIES = {Springer Monographs in Mathematics},
      NOTE = {A view from variational analysis},
 PUBLISHER = {Springer, Dordrecht},
      YEAR = {2009},
     PAGES = {xii+375},
      ISBN = {978-0-387-87820-1},
   MRCLASS = {49J53 (26B10 47J30 47N10 49K40 58C15 90C31)},
  MRNUMBER = {2515104},
MRREVIEWER = {Stephen\ M.\ Robinson},
       DOI = {10.1007/978-0-387-87821-8},
       URL = {https://doi.org/10.1007/978-0-387-87821-8},
}

@article {DrusvyatskiyLewis,
    	AUTHOR = {Drusvyatskiy, D. and Lewis, A. S.},
     	TITLE = {Tilt stability, uniform quadratic growth, and strong metric regularity of the subdifferential},
   	JOURNAL = {SIAM J. Optim.},
  	FJOURNAL = {SIAM Journal on Optimization},
    	VOLUME = {23},
      	YEAR = {2013},
    	NUMBER = {1},
     	PAGES = {256--267},
      	ISSN = {1052-6234,1095-7189},
   	MRCLASS = {49J53 (49J52 49K40 90C31)},
  	MRNUMBER = {3033107},
	MRREVIEWER = {Shawn\ Xianfu\ Wang},
       	DOI = {10.1137/120876551},
       	URL = {https://doi.org/10.1137/120876551},
}

@article {DMN,
    	AUTHOR = {Drusvyatskiy, D. and Mordukhovich, B. S. and Nghia, T. T. A.},
     	TITLE = {Second-order growth, tilt stability, and metric regularity of the subdifferential},
   	JOURNAL = {J. Convex Anal.},
  	FJOURNAL = {Journal of Convex Analysis},
    	VOLUME = {21},
      	YEAR = {2014},
    	NUMBER = {4},
     	PAGES = {1165--1192},
      	ISSN = {0944-6532,2363-6394},
   	MRCLASS = {49J52 (49J53 49K40 90C31)},
  	MRNUMBER = {3331214},
	MRREVIEWER = {Francisco\ J.\ Arag\'on Artacho},
}

@article {HangSarabi,
    	AUTHOR = {Hang, Nguyen T. V. and Sarabi, Ebrahim},
     	TITLE = {Smoothness of subgradient mappings and its applications in parametric optimization},
   	JOURNAL = {Set-Valued Var. Anal.},
  	FJOURNAL = {Set-Valued and Variational Analysis. Theory and Applications},
    	VOLUME = {33},
      	YEAR = {2025},
    	NUMBER = {4},
     	PAGES = {Paper No. 41, 40},
      	ISSN = {1877-0533,1877-0541},
   	MRCLASS = {90C31 (49J52 49J53 65K15)},
  	MRNUMBER = {4963701},
	MRREVIEWER = {Juan\ Parra},
       	DOI = {10.1007/s11228-025-00777-z},
       	URL = {https://doi.org/10.1007/s11228-025-00777-z},
}

@unpublished{MaSun,
    	AUTHOR = {Ma, J. and Sun, D.},
	    	TITLE = {The {Aubin} property for generalized equations over {$C^2$}-cone reducible sets},
    	NOTE = {arXiv preprint, arXiv:2509.14194v5},
    	YEAR = {2026},
}

@article{AndreaniHaeserRamosSilva,
        AUTHOR = {Andreani, Roberto and Haeser, Gabriel and Ramos, Alberto and Silva, Paulo J. S.},
         TITLE = {A Second-Order Sequential Optimality Condition Associated to the Convergence of Optimization Algorithms},
       JOURNAL = {IMA J. Numer. Anal.},
        VOLUME = {37},
          YEAR = {2017},
        NUMBER = {4},
         PAGES = {1902--1929},
           DOI = {10.1093/imanum/drw064},
           URL = {https://doi.org/10.1093/imanum/drw064},
}

@article{FukudaHaeserMito,
        AUTHOR = {Fukuda, Ellen H. and Haeser, Gabriel and Mito, Leonardo M.},
         TITLE = {On the Weak Second-Order Optimality Condition for Nonlinear Semidefinite and Second-Order Cone Programming},
       JOURNAL = {Set-Valued Var. Anal.},
        VOLUME = {31},
          YEAR = {2023},
        NUMBER = {2},
         PAGES = {Paper No. 15, 28},
           DOI = {10.1007/s11228-023-00676-1},
           URL = {https://doi.org/10.1007/s11228-023-00676-1},
}

@article{DavisDrusvyatskiy,
        AUTHOR = {Davis, Damek and Drusvyatskiy, Dmitriy},
         TITLE = {Proximal Methods Avoid Active Strict Saddles of Weakly Convex Functions},
       JOURNAL = {Found. Comput. Math.},
        VOLUME = {22},
          YEAR = {2022},
        NUMBER = {2},
         PAGES = {561--606},
           DOI = {10.1007/s10208-021-09516-w},
           URL = {https://doi.org/10.1007/s10208-021-09516-w},
}

@article{BianchiHachemSchechtman,
        AUTHOR = {Bianchi, Pascal and Hachem, Walid and Schechtman, Sholom},
         TITLE = {Stochastic Subgradient Descent Escapes Active Strict Saddles on Weakly Convex Functions},
       JOURNAL = {Math. Oper. Res.},
        VOLUME = {49},
          YEAR = {2024},
        NUMBER = {3},
         PAGES = {1761--1790},
           DOI = {10.1287/moor.2021.0194},
           URL = {https://doi.org/10.1287/moor.2021.0194},
}

@article {PolRoc96,
    AUTHOR = {Poliquin, R. A. and Rockafellar, R. T.},
     TITLE = {Prox-regular functions in variational analysis},
   JOURNAL = {Trans. Amer. Math. Soc.},
  FJOURNAL = {Transactions of the American Mathematical Society},
    VOLUME = {348},
      YEAR = {1996},
    NUMBER = {5},
     PAGES = {1805--1838},
      ISSN = {0002-9947,1088-6850},
   MRCLASS = {49J52 (26B25 54C60 58C06)},
  MRNUMBER = {1333397},
MRREVIEWER = {Philip\ D.\ Loewen},
       DOI = {10.1090/S0002-9947-96-01544-9},
       URL = {https://doi.org/10.1090/S0002-9947-96-01544-9},
}

@book {RockafellarConvexAnalysis,
        AUTHOR = {Rockafellar, R. Tyrrell},
         TITLE = {Convex analysis},
        SERIES = {Princeton Mathematical Series},
        VOLUME = {28},
     PUBLISHER = {Princeton University Press, Princeton, NJ},
          YEAR = {1970},
         PAGES = {xviii+451},
          ISBN = {0-691-08069-0},
      MRNUMBER = {0274683},
}

@book {FenchelConvexCones,
        AUTHOR = {Fenchel, Werner},
         TITLE = {Convex cones, sets, and functions},
     PUBLISHER = {Department of Mathematics, Princeton University},
       ADDRESS = {Princeton, NJ},
          YEAR = {1953},
          NOTE = {Notes by D. W. Blackett of lectures given at Princeton
                  University in the spring term of 1951},
           URL = {https://convexoptimization.com/TOOLS/Fenchel1951.pdf},
}

@book {RockafellarWets,
    	AUTHOR = {Rockafellar, R. Tyrrell and Wets, Roger J.-B.},
     	TITLE = {Variational analysis},
    	SERIES = {Grundlehren der mathematischen Wissenschaften [Fundamental Principles of Mathematical Sciences]},
    	VOLUME = {317},
 	PUBLISHER = {Springer-Verlag, Berlin},
      	YEAR = {1998},
     	PAGES = {xiv+733},
      	ISBN = {3-540-62772-3},
   	MRCLASS = {49-02 (46N10 47N10 49J52 49K40 90C30)},
  	MRNUMBER = {1491362},
	MRREVIEWER = {Francis\ H.\ Clarke},
       	DOI = {10.1007/978-3-642-02431-3},
       	URL = {https://doi.org/10.1007/978-3-642-02431-3},
}

@book {vanDenDries,
    	AUTHOR = {van den Dries, Lou},
     	TITLE = {Tame topology and o-minimal structures},
    	SERIES = {London Mathematical Society Lecture Note Series},
    	VOLUME = {248},
 	PUBLISHER = {Cambridge University Press, Cambridge},
      	YEAR = {1998},
     	PAGES = {x+180},
      	ISBN = {0-521-59838-9},
   	MRCLASS = {03-02 (03C50 03C60 14P10 52-02 54-02 55-02 57-02)},
  	MRNUMBER = {1633348},
	MRREVIEWER = {O.\ V.\ Belegradek},
       	DOI = {10.1017/CBO9780511525919},
        URL = {https://doi.org/10.1017/CBO9780511525919},
}

@book {Hatcher,
        AUTHOR = {Hatcher, Allen},
         TITLE = {Algebraic topology},
     PUBLISHER = {Cambridge University Press, Cambridge},
          YEAR = {2002},
         PAGES = {xii+544},
          ISBN = {0-521-79540-0; 0-521-79160-X},
           URL = {https://pi.math.cornell.edu/~hatcher/AT/AT.pdf},
}

@book {Munkres,
        AUTHOR = {Munkres, James R.},
         TITLE = {Topology},
       EDITION = {Second},
     PUBLISHER = {Prentice Hall, Inc., Upper Saddle River, NJ},
          YEAR = {2000},
         PAGES = {xvi+537},
          ISBN = {0-13-181629-2},
}

@article {Dutertre2020,
        AUTHOR = {Dutertre, Nicolas},
         TITLE = {On the topology of non-isolated real singularities},
       JOURNAL = {Journal of Singularities},
        VOLUME = {22},
          YEAR = {2020},
         PAGES = {159--179},
           DOI = {10.5427/jsing.2020.22j},
           URL = {https://doi.org/10.5427/jsing.2020.22j},
}

@article {Escribano2002,
        AUTHOR = {Escribano, Jes{\'u}s},
         TITLE = {Approximation theorems in o-minimal structures},
       JOURNAL = {Illinois J. Math.},
        VOLUME = {46},
          YEAR = {2002},
        NUMBER = {1},
         PAGES = {111--128},
          ISSN = {0019-2082},
           DOI = {10.1215/ijm/1258136143},
           URL = {https://doi.org/10.1215/ijm/1258136143},
}

@article {Tanabe2023,
        AUTHOR = {Tanabe, Masato},
         TITLE = {Canonical stratification of definable {L}ie groupoids},
       JOURNAL = {Journal of Singularities},
        VOLUME = {26},
          YEAR = {2023},
         PAGES = {63--75},
          ISSN = {1949-2006},
           DOI = {10.5427/jsing.2023.26d},
           URL = {https://doi.org/10.5427/jsing.2023.26d},
}

@article {Mather2012,
        AUTHOR = {Mather, John N.},
         TITLE = {Notes on topological stability},
       JOURNAL = {Bull. Amer. Math. Soc. (N.S.)},
        VOLUME = {49},
          YEAR = {2012},
        NUMBER = {4},
         PAGES = {475--506},
          ISSN = {0273-0979},
           DOI = {10.1090/S0273-0979-2012-01383-6},
           URL = {https://doi.org/10.1090/S0273-0979-2012-01383-6},
}

@book {Hirsch,
        AUTHOR = {Hirsch, Morris W.},
         TITLE = {Differential topology},
        SERIES = {Graduate Texts in Mathematics},
        VOLUME = {33},
     PUBLISHER = {Springer-Verlag, New York},
          YEAR = {1976},
         PAGES = {x+221},
          ISBN = {0-387-90148-5},
           DOI = {10.1007/978-1-4684-9449-5},
           URL = {https://doi.org/10.1007/978-1-4684-9449-5},
}

@book {MargalefOuterelo,
        AUTHOR = {Margalef-Roig, Juan and Outerelo Dom{\'i}nguez, Enrique},
         TITLE = {Differential topology},
        SERIES = {North-Holland Mathematics Studies},
        VOLUME = {173},
     PUBLISHER = {North-Holland Publishing Co., Amsterdam},
          YEAR = {1992},
         PAGES = {xvi+604},
          ISBN = {0-444-88434-3},
      MRNUMBER = {1173211},
}

@article {Jubin2023,
        AUTHOR = {Jubin, Beno{\^\i}t},
         TITLE = {Intrinsic volumes of sublevel sets},
       JOURNAL = {Ann. Fac. Sci. Toulouse Math. (6)},
        VOLUME = {32},
          YEAR = {2023},
        NUMBER = {5},
         PAGES = {911--938},
           DOI = {10.5802/afst.1758},
           URL = {https://doi.org/10.5802/afst.1758},
}

@article {GalashinKarpLam2022,
        AUTHOR = {Galashin, Pavel and Karp, Steven N. and Lam, Thomas},
         TITLE = {Regularity theorem for totally nonnegative flag varieties},
       JOURNAL = {J. Amer. Math. Soc.},
        VOLUME = {35},
          YEAR = {2022},
        NUMBER = {2},
         PAGES = {513--579},
           DOI = {10.1090/jams/983},
           URL = {https://doi.org/10.1090/jams/983},
}

@article {BenkoMehlitz2023,
        AUTHOR = {Benko, Mat{\'u}\v{s} and Mehlitz, Patrick},
         TITLE = {Why Second-Order Sufficient Conditions Are, in a Way,
                  Easy---or---Revisiting the Calculus for Second
                  Subderivatives},
       JOURNAL = {J. Convex Anal.},
        VOLUME = {30},
          YEAR = {2023},
        NUMBER = {2},
         PAGES = {541--589},
           URL = {https://www.heldermann.de/JCA/JCA30/JCA302/jca30031.htm},
}

@book {HornJohnson2013,
        AUTHOR = {Horn, Roger A. and Johnson, Charles R.},
         TITLE = {Matrix Analysis},
       EDITION = {Second},
     PUBLISHER = {Cambridge University Press},
       ADDRESS = {Cambridge},
          YEAR = {2013},
         PAGES = {xviii+643},
          ISBN = {978-0-521-83940-2; 978-0-521-54823-6},
           DOI = {10.1017/CBO9781139020411},
           URL = {https://doi.org/10.1017/CBO9781139020411},
}

@article {RockafellarHiddenConvexity,
        AUTHOR = {Rockafellar, R. Tyrrell},
         TITLE = {Augmented Lagrangians and hidden convexity in sufficient
                  conditions for local optimality},
       JOURNAL = {Math. Program.},
        VOLUME = {198},
          YEAR = {2023},
        NUMBER = {1},
         PAGES = {159--194},
           DOI = {10.1007/s10107-022-01768-w},
           URL = {https://doi.org/10.1007/s10107-022-01768-w},
}

@article {Haraux2005,
        AUTHOR = {Haraux, Alain},
         TITLE = {Positively homogeneous functions and the
                  {{\L}ojasiewicz} gradient inequality},
       JOURNAL = {Ann. Polon. Math.},
        VOLUME = {87},
          YEAR = {2005},
         PAGES = {165--174},
           DOI = {10.4064/ap87-0-13},
           URL = {https://doi.org/10.4064/ap87-0-13},
}

@book{absil2008manifold,
    AUTHOR = {Absil, P.-A. and Mahony, Robert and Sepulchre, Rodolphe},
     TITLE = {Optimization algorithms on matrix manifolds},
 PUBLISHER = {Princeton University Press, Princeton, NJ},
      YEAR = {2008},
      ISBN = {978-0-691-13298-3},
       URL = {https://press.princeton.edu/books/ebook/9781400830244/optimization-algorithms-on-matrix-manifolds-pdf},
}

\end{document}